\documentclass{ws-idaqp}

\usepackage[latin1]{inputenc}

\usepackage{amsmath}
\usepackage{amssymb}
\usepackage[mathscr]{eucal}

    \usepackage{theorem}

\newenvironment{env}[2]{\begin{#1}#2\end{#1}}{}
    \newcommand{\beq}[1]{\begin{env}{equation}{#1}}
    \newcommand{\beqn}[1]{\begin{env}{equation*}{#1}}
    \newcommand{\bal}[1]{\begin{env}{align}{#1}}
    \newcommand{\baln}[1]{\begin{env}{align*}{#1}}
    \newcommand{\bga}[1]{\begin{env}{gather}{#1}}
    \newcommand{\bgan}[1]{\begin{env}{gather*}{#1}}
    \newcommand{\bflal}[1]{\begin{env}{flalign}{#1}}
    \newcommand{\bflaln}[1]{\begin{env}{flalign*}{#1}}
    \newcommand{\bmu}[1]{\begin{env}{multline}{#1}}
    \newcommand{\bmun}[1]{\begin{env}{multline*}{#1}}
    \newcommand{\bsp}[1]{\begin{env}{split}{#1}}

    \newcommand{\eeq}{\end{env}}
    \newcommand{\eeqn}{\end{env}}
    \newcommand{\eal}{\end{env}}
    \newcommand{\ealn}{\end{env}}
    \newcommand{\ega}{\end{env}}
    \newcommand{\egan}{\end{env}}
    \newcommand{\eflal}{\end{env}}
    \newcommand{\eflaln}{\end{env}}
    \newcommand{\emu}{\end{env}}
    \newcommand{\emun}{\end{env}}
    \newcommand{\esp}{\end{env}}

\newcommand{\lf}{\vspace{2ex}}

\renewcommand{\bf}[1]{\textbf{#1}}
\renewcommand{\it}[1]{\textit{#1}}
\renewcommand{\sc}[1]{\textsc{#1}}
\renewcommand{\sf}[1]{\textsf{#1}}
\renewcommand{\sl}[1]{\textsl{#1}}
\renewcommand{\tt}[1]{\texttt{#1}}
\newcommand{\hl}[1]{\bf{\it{#1}}}
\newcommand{\mrm}[1]{\mathrm{#1}}
\newcommand{\mbf}[1]{\mathbf{#1}}
\newcommand{\msf}[1]{\text{\small$\sf{#1}$}}

\newcommand{\cmc}[1]{\mathcal{#1}}
\newcommand{\eus}[1]{\mathscr{#1}}
\newcommand{\euf}[1]{\mathfrak{#1}}
\newcommand{\bb}[1]{\mathbb{#1}}

\newcommand{\mscriptsize}[1]{{\setlength{\arraycolsep}{.3ex}\text{\scriptsize$#1$}}}

\newcommand{\nbd}[1]{$#1$\nobreakdash--}
\newcommand{\ol}[1]{\overline{#1}}
\newcommand{\ul}[1]{\underline{#1}}

\newcommand{\wh}[1]{\widehat{#1}}
\newcommand{\ve}{\varepsilon}
\newcommand{\vt}{\vartheta}
\newcommand{\vk}{\varkappa}
\newcommand{\vp}{\varphi}
\newcommand{\om}{\omega}

\newcommand{\abs}[1]{\left\lvert#1\right\rvert}
\newcommand{\norm}[1]{\left\lVert#1\right\rVert}

\newcommand{\bfam}[1]{\bigl(#1\bigr)}
\newcommand{\Bfam}[1]{\Bigl(#1\Bigr)}
\newcommand{\AB}[1]{\langle#1\rangle}
\newcommand{\bAB}[1]{\bigl\langle#1\bigr\rangle}
\newcommand{\BAB}[1]{\Bigl\langle#1\Bigr\rangle}
\newcommand{\CB}[1]{\{#1\}}
\newcommand{\bCB}[1]{\bigl\{#1\bigr\}}
\newcommand{\BCB}[1]{\Bigl\{#1\Bigr\}}
\newcommand{\SB}[1]{[#1]}

\newcommand{\RO}[1]{[#1)}

\newcommand{\Matrix}[1]{\begin{pmatrix}#1\end{pmatrix}}

\newcommand{\sMatrix}[1]{\mscriptsize{\Matrix{#1}}}

\newcommand{\sbar}[1]{\:\bar{#1}\:}

\newcommand{\set}[2][]{
    \ifthenelse{\equal{#1}{}}{
        \CB{#2}}{
        \CB{#1~|~#2}}}
\newcommand{\bset}[2][]{
    \ifthenelse{\equal{#1}{}}{
        \bCB{#2}}{
        \bCB{#1~|~#2}}}
\newcommand{\Bset}[2][]{
    \ifthenelse{\equal{#1}{}}{
        \BCB{#2}}{
        \BCB{#1~\big|~#2}}}
\newcommand{\zero}{\CB{0}}

\DeclareMathOperator{\ls}{\normalfont\msf{span}}
\DeclareMathOperator{\cls}{\ol{\ls}}
\DeclareMathOperator*{\limind}{lim\,ind}

\DeclareMathOperator*{\coplus}{\ol{\bigoplus}}

\DeclareMathOperator{\id}{\normalfont\msf{id}}

\renewcommand{\dim}{\operatorname{\msf{dim}}}

\newcommand{\C}{\bb{C}}

\newcommand{\bI}{\bb{I}}
\newcommand{\bJ}{\bb{J}}

\newcommand{\N}{\bb{N}}

\newcommand{\R}{\bb{R}}

\newcommand{\cB}{\cmc{B}}

\newcommand{\cF}{\cmc{F}}

\newcommand{\sB}{\eus{B}}

\newcommand{\sS}{\eus{S}}

\newcommand{\sU}{\eus{U}}

\newcommand{\er}{\euf{r}}
\newcommand{\es}{\euf{s}}
\newcommand{\et}{\euf{t}}

\newcommand{\eB}{\euf{B}}

\newcommand{\eH}{\euf{H}}

\newcommand{\eK}{\euf{K}}
\newcommand{\eL}{\euf{L}}

\newcommand{\eU}{\euf{U}}

\newcommand{\U}{\mbf{1}}

\newcommand{\G}{\Gamma}

\newcommand{\DG}{{\mrm{I}\hspace{-0.3ex}\G}}

\newcommand{\I}{{I\!\!\!\;I}}

\newcommand{\s}{\text{\scriptsize$\sS$}}

\newcommand{\bplus}{\text{\,\footnotesize$\boxplus$\,}}

    \numberwithin{equation}{section}
    \renewcommand{\theequation}{\thesection.\arabic{equation}}

        \newcommand{\mnname}{Mathematical note.}
        \newcommand{\enname}{End of the note.}
        \newcommand{\definame}{Definition.}
        \newcommand{\propname}{Proposition.}
        \newcommand{\lemname}{Lemma.}
        \newcommand{\exname}{Example.}
        \newcommand{\exername}{Exercise.}
        \newcommand{\remname}{Remark.}
        \newcommand{\obname}{Observation.}
        \newcommand{\thmname}{Theorem.}
        \newcommand{\corname}{Corollary.}

    \theoremheaderfont{\normalfont\bfseries}
    \theoremstyle{change}
        \theorembodyfont{\rmfamily}
            \newtheorem{emp}{}[section]
                \newcommand{\bemp}[1][]{
                    \begin{emp}\hskip-\labelsep\bf{#1}\hskip\labelsep}
                \newcommand{\eemp}{\end{emp}}
\newtheorem{itemp}[emp]{}
                \newcommand{\bitemp}[1][]{
                    \begin{itemp}\hskip-\labelsep\bf{#1}\hskip\labelsep\normalfont\itshape}
                \newcommand{\eitemp}{\end{itemp}}
            \newtheorem{mn}[emp]{\mnname}
                \newcommand{\bmn}{\begin{mn}~\begin{quotation}\renewcommand{\baselinestretch}{1}\small\noindent\ignorespaces}
                \newcommand{\emn}{\end{quotation}\hfill\bf{\enname}\end{mn}}
            \newtheorem{ex}[emp]{\exname}
                \newcommand{\bex}{\begin{ex}}
                \newcommand{\eex}{\end{ex}}
            \newtheorem{exer}[emp]{\exername}
                \newcommand{\bexer}{\begin{exer}}
                \newcommand{\eexer}{\end{exer}}
            \newtheorem{defi}[emp]{\definame}
                \newcommand{\bdefi}{\begin{defi}}
                \newcommand{\edefi}{\end{defi}}
            \newtheorem{rem}[emp]{\remname}
                \newcommand{\brem}{\begin{rem}}
                \newcommand{\erem}{\end{rem}}
            \newtheorem{ob}[emp]{\obname}
                \newcommand{\bob}{\begin{ob}}
                \newcommand{\eob}{\end{ob}}
        \theorembodyfont{\normalfont\itshape}
            \newtheorem{thm}[emp]{\thmname}
                \newcommand{\bthm}{\begin{thm}}
                \newcommand{\ethm}{\end{thm}}
            \newtheorem{prop}[emp]{\propname}
                \newcommand{\bprop}{\begin{prop}}
                \newcommand{\eprop}{\end{prop}}
            \newtheorem{cor}[emp]{\corname}
                \newcommand{\bcor}{\begin{cor}}
                \newcommand{\ecor}{\end{cor}}
            \newtheorem{lem}[emp]{\lemname}
                \newcommand{\blem}{\begin{lem}}
                \newcommand{\elem}{\end{lem}}
\newenvironment{empn}[1]{\lf\noindent\bf{#1}\ignorespaces\hskip\labelsep}{\lf}
		\newcommand{\bempn}[1]{\begin{empn}{#1}}
		\newcommand{\eempn}{\end{empn}}
		\newcommand{\bitempn}[1]{\begin{empn}{#1}\normalfont\itshape}
		\newcommand{\eitempn}{\end{empn}}
                \newcommand{\bmnn}{\begin{empn}{\mnname}~\begin{quotation}\renewcommand{\baselinestretch}{1}\small\noindent\ignorespaces}
                \newcommand{\emnn}{\end{quotation}\hfill\bf{\enname}\end{empn}}
		\newcommand{\bexn}{\begin{empn}{\exname}}
		\newcommand{\eexn}{\end{empn}}
		\newcommand{\bexern}{\begin{empn}{\exername}}
		\newcommand{\eexern}{\end{empn}}
		\newcommand{\bdefin}{\begin{empn}{\definame}}
		\newcommand{\edefin}{\end{empn}}
		\newcommand{\bremn}{\begin{empn}{\remname}}
		\newcommand{\eremn}{\end{empn}}
		\newcommand{\bobn}{\begin{empn}{\obname}}
		\newcommand{\eobn}{\end{empn}}

\usepackage{float}
\usepackage{epic}
\usepackage{eepic}
\usepackage[matrix,arrow]{xy}

\renewcommand{\msf}[1]{\mathsf{#1}}

\begin{document}



\renewcommand{\thefootnote}{{(\alph{footnote})}}



\title{The Index of (White) Noises and
\\
their Product Systems\thanks{This work is supported by a PPP-project by DAAD and DST.}}

\author{Michael Skeide}

\address{Dipartimento S.E.G.e S., Università degli Studei del Molise, Via de Sancits\\86100 Campobasso, Italy\\e-mail:\tt{skeide@unimol.it}}


\maketitle


\begin{abstract}
\noindent
Almost every article about \it{Arveson systems} (that is, product systems of Hilbert spaces) starts by recalling their basic classification assigning to every Arveson system a type and an index. So it is natural to ask in how far an analogue classification can be proposed also for product systems of Hilbert modules. However, while the definition of type is plain, there are obstacles for the definition of index. But all obstacles can be removed when restricting to the category which we introduce here as spatial product systems and that matches the usual definition of spatial in the case of Arveson systems. This is not really a loss because the definition of index for nonspatial Arveson systems is rather formal and does not reflect the information the index carries for spatial Arveson systems.

\nbd{E_0}semigroups give rise to product systems. Our definition of spatial product system, namely, existence of a unital unit that is central, matches Powers' definition of spatial in the sense that the \nbd{E_0}semigroup from which the product system is derived admits a semigroup of intertwining isometries. We show that every spatial product system contains a unique maximal completely spatial subsystem (generated by all units) that is isomorphic to a product system of time ordered Fock modules. (There exist nonspatial product system that are generated by their units. Consequently, these cannot be Fock modules.) The index of a spatial product system we define as the (unique) Hilbert bimodule that determins the Fock module. In order to show that the index merits the name index we provide a product of product systems under which the index is additive (direct sum). While for Arveson systems there is the tensor product, for general product systems the tensor product does not make sense as a product system. Even for Arveson systems our product is, in general, only a subsystem of the tensor product. Moreover, its construction depends explicitly on the choice of the central reference units of its factors.

Spatiality of a product system means that it may be derived from an \nbd{E_0}semigroup with an invariant vector expectation, that is, from a noise. We extend our product of spatial product systems to a product of noises and study its properties.

Finally, we apply our techniques to show the module analogue of Fowler's result that free flows are comletely spatial, and we compute their indices.
\end{abstract}




\section{Introduction} \label{intro}

{\parskip0.5ex plus 0.5ex minus 0.5ex
\noindent{\bf{Arveson systems and \nbd{E_0}semigroups.~}}
In a series of papers (see in particular [\refcite{Arv89,Arv89a}]) Arveson worked out a close relationship between \it{$E_0$--semi\-groups} (semigroups of unital endomorphisms of a $C^*$--al\-ge\-bra) on $\sB(H)$ and \it{product systems of Hilbert spaces} (Arveson systems for short). More precisely, he discovered (we do not repeat the precise technical conditions) a one-to-one correspondence between \nbd{E_0}semigroups (up to \it{cocycle conjugacy}) and Arveson system (up to isomorphism).

An Arveson system is a family of $\eH^\otimes=\bfam{\eH_t}_{t\in\R_+}$ of Hilbert spaces $\eH_t$ fulfilling
\beqn{
\eH_{s+t}
~=~
\eH_s{\sbar{\otimes}}\eH_t
}\eeqn
in an associative way. (Actually, there are also some measurability requirements and the $\eH_t$ should be infinite-dimensional and separable for $t>0$, but we do not speak about this.) Arveson introduced also the concept of \it{units}, i.e.\ a family $u^\otimes=\bfam{u_t}_{t\in\R_+}$ of nonzero elements $u_t\in\eH_t$ fulfilling
\beqn{
u_{s+t}
~=~
u_s\otimes u_t.
}\eeqn
(Once again there is a measurability condition which we ignore.) The construction of the Arveson system associated with an \nbd{E_0}semigroup is plain and if an Arveson system has a unit, then it is also easy to construct an \nbd{E_0}semigroup that has associated with it the Arveson system we started with. (We discuss the more general version for Hilbert modules detailed in Section \ref{wnspsSEC}.) The construction of an \nbd{E_0}semigroup from a unitless Arveson system in [\refcite{Arv89a}] is among the most difficult results about Arveson systems.

The simplest example of an Arveson system is the family $\G^\otimes(K)=\bfam{\G_t(K)}_{t\in\R_+}$ of symmetric Fock spaces $\G_t(K)=\G(L^2(\SB{0,t},K))$ with isomorphism
\beqn{
\G_{s+t}(K)
~=~
\G(L^2(\SB{t,t+s},K)){\sbar{\otimes}}\G_t(K)
~=~
\G_s(K){\sbar{\otimes}}\G_t(K).
}\eeqn
The (measurable) units are precisely those given by $u_t=e^{tc}\psi(\I_{\SB{0,t}}f)$ where $c\in\C$ and $\psi(\I_{\SB{0,t}}f)$ is the exponential vector to the funcion $\I_{\SB{0,t}}f$ with $f\in K$. (This product system belongs, for instance, to the \nbd{E_0}semigroup on $\sB\bfam{\G(L^2(\R_+,K))}$ induced by the time shift on $\G(L^2(\R_+,K))$.) Arveson showed that all product systems of Hilbert spaces which are spanned linearly by tensor products of their units (so-called \it{type I systems}) are isomorphic to some $\G^\otimes(K)$ for a suitable Hilbert space $K$.

Arveson systems which have a unit are called \it{spatial}. In general, any Arveson system contains a maximal type I subsystem (namely, that which is generated be the units). The \it{Arveson index} of an Arveson system is the dimension of the the Hilbert spaces $K$ for the maximal type I subsystem. It is put to $\infty$ by hand, if the subsystem is $\CB{0}$, i.e.\ for nonspatial Arveson systems. (We explain in Section \ref{indnonspat} why we think the index should be defined only for spatial Arveson systems.) The index of an \nbd{E_0}semigroup is that of its associated Arveson system.

Be it among \nbd{E_0}semigroups on $\sB(H)$ be it among Arveson systems there is a natural operation, the tensor product, making out of two of them a new one. Obviously, $\G^\otimes(K){\sbar{\otimes}}\G^\otimes(K')=\G^\otimes(K\oplus K')$ so that the index is additive under tensor product and, thus, indeed merits to be named \it{index}.

The index is a complete isomorphism invariant for type I (or \it{completely spatial}) Arveson systems. For other spatial Arveson systems this is not so (and for nonspatial systems the index has no good meaning). Appart from Powers' examples (see, for instance, [\refcite{Pow87}]) it was Tsirelson who, using probabilistic ideas, provided us with larger classes of examples (see [\refcite{Tsi00p1}] for spatial Arveson systems and [\refcite{Tsi00p2}] for nonspatial examples).\footnote{Liebscher [\refcite{Lie00p1}] started to provide us with more (still incomplete) isomorphism invariants for spatial Arveson systems based on substantial extensions of the ideas of [\refcite{Tsi00p1}], while Bhat and Srinivasan [\refcite{BhSr05}] initiated a systematic study of nonspatial Arveson system using a more funcional analytic approach to the ideas of [\refcite{Tsi00p2}]. In [\refcite{Pow03}] Powers reduced the study of spatial \nbd{E_0}semigroups and, therefore, of spatial Arveson systems to the study of so-called \it{CP-flows}.}

\lf\noindent
\bf{Product systems of Hilbert modules.~}
General product systems of Hilbert \nbd{\cB}\nbd{\cB}modules occured in Bhat and Skeide [\refcite{BhSk00}] in dilation theory of CP-semigroups on a unital \nbd{C^*}algebra $\cB$ dilating the CP-semigroup to an \nbd{E_0}semigroup on $\sB^a(E)$ for some Hilbert \nbd{\cB}module $E$.\footnote{In the case of von Neumann algebras Muhly and Solel [\refcite{MuSo02}] have constructed the same dilation with the help of a product system of von Neumann \nbd{\cB'}\nbd{\cB'}modules, where $\cB'$ is the commutant of $\cB$. The duality between the two approaches (in fact, the construction of a \it{commutant} of von Neumann bimodules as introduced in Skeide [\refcite{Ske03c}] and also in Muhly and Solel [\refcite{MuSo04}]) is explained in Skeide [\refcite{Ske03c}].} The technical definition of product systems (continuity or measurability conditions) will depend on the purpose.\footnote{See for instance Skeide [\refcite{Ske03b}] or Hirshberg [\refcite{Hir04,Hir05}].} In fact, we prefer to investigate always also the algebraic case, in order not to exclude interesting product systems (that exist!) from the discussion just because we are not (yet) able to show that they fulfill certain (possibly premature) technical conditions.

It is the goal of these notes to see in how far it is possible to obtain the basic results about spatial Arveson system also for spatial product systems of Hilbert modules. This requires, in the first place, to single out the correct notion of \it{spatial} product system. In order to define an index we must specify \it{completely spatial} product systems and see whether they are isomorphic to some sort of Fock module that substitutes the symmetric Fock space. Finally, we must find a \it{product} of spatial product systems that substitutes the tensor product of (spatial) Arveson systems, because the tensor product of product systems does not make sense, in general, for Hilbert modules.

For that goal we follow the theory of product systems as far as possible in analogy with the theory of Arveson systems. We consider all product systems as derived from \nbd{E_0}semigroups by the construction from Skeide [\refcite{Ske02}] (a straightforward generalization to Hilbert modules of Bhat's approach to Arveson systems in [\refcite{Bha96}]). It is one of the big open questions, whether every product system can be obtained in that way, but if the product system has a unit, like our spatial product systems, then it is true; see Section \ref{wnspsSEC} for a detailed explanation.

Powers [\refcite{Pow87}] calls an \nbd{E_0}semigroup on $\sB(H)$ spatial, if it admits an intertwining semigroup of isometries. It is easy to give examples\footnote{See Skeide [\refcite{Ske03b}]} that in the case of Hilbert \nbd{\cB}modules mere existence of a unit in a product system is not sufficient to achieve this. The unit must be \it{central}, i.e.\ its members must commute with the elements of $\cB$, --- and this is our definition (Definition \ref{spatdef}) of spatial product systems.

The \it{time ordered Fock module} over a Hilbert \nbd{\cB}\nbd{\cB}module $F$ (introduced in Bhat and Skeide [\refcite{BhSk00}] and studied in detail in Liebscher and Skeide [\refcite{LiSk01}]) is the analogue of the symmetric Fock space. A time ordered product system has a central unital unit (the vacuum), so that it is spatial, and it is generated by its units, so that it is even \it{completely spatial}. Using results from [\refcite{LiSk01,BBLS04}] we show (Theorem \ref{msssthm}) that (in analogy with Arveson systems) every completely spatial product system is isomorphic to a time ordered system. Also here it is easy to give counter examples which show that existence of a cenral unit may not be dropped from the definition of completely spatial product system. (A product system that is generated by its units is a strongly dense subsystem of a time ordered product system of von Neumann modules over the enveloping von Neumann algebra of $\cB$. In fact, it is the main result of Barreto, Bhat, Liebscher and Skeide [\refcite{BBLS04}] that every product system of von Neumann modules that has a unit is spatial. This result is equivalent to the results by Christensen and Evans [\refcite{ChrEv79}] on the form of the generator of uniformly continuous normal CP-semigroups on a von Neumann algebra.)

Spatial product systems contain a unique maximal completely spatial subsystem (Corollary \ref{maxspat}). By Theorem \ref{msssthm} the completely spatial subsystem is isomorphic to a time ordered system and this time ordered system is determined by the unique Hilbert bimodule $F$ that plays the role of the Hilbert space $K$ for type I Arveson systems. As $F$ is no longer determined by a simple dimension we stick to the whole space $F$ as \it{index} of the spatial product system (Definition \ref{indexdef}). We construct a product of spatial product systems (Definition \ref{spatproddef}) under which the indices of the factors add up as direct sums (Theorem \ref{indsumthm}). In the case of spatial Arveson systems our product is a subsystem of the tensor product that may but need not coincide with the tensor product. (In fact, by Theorem \ref{prodpropthm} our product is generated by its factors, while the tensor product of Arveson systems need not be.)\footnote{In Section \ref{powersprod} we report an example of Powers (after publication of these notes) where our product naturally occurs and is not the tensor product.}

All spatial product systems can be derived from an \nbd{E_0}semigroup with an invariant conditional vector expectation, that is from \it{noises}. We extend our product of spatial product systems to a product of noises in such a way that the associated spatial product system is the product of the spatial product systems associated with the factors (Theorem \ref{WNprodthm}). We add the result (Theorem \ref{scprodthm}) that the product of noises preserves strong continuity in time.\footnote{Together with the results from Skeide [\refcite{Ske03b}] this implies also that the product of continuous (in the sense of [\refcite{Ske03b}]) spatial product systems is continuous, that is our product is compatible with technical definitions of product systems.}

As a concrete example we show that the time shift semigroup on the \it{full Fock module} is completely spatial and we calculate its index (Theorem \ref{Fowmodthm}). Motivated by the fact that these {free flows} have sitting inside also a \it{free product system}, we suggest a couple of natural questions arround free product systems.

\lf\noindent
\bf{Contents.~}
These notes are organized as follows. Sections \ref{psSEC} and \ref{wnspsSEC} are quasi completely repetitive. In Section \ref{psSEC} we recall basic definitions concerning product systems like units and morphisms. Then we define spatial product systems and their morphisms. These definitions should suffice to understand the algebraic construction of the product of spatial product systems in Section \ref{prodSEC}. In Section \ref{wnspsSEC} we repeat the relation between product systems and \nbd{E_0}semigroups and point at some specific properties in the case of spatial product systems and noises. Sections \ref{psSEC}, \ref{prodSEC} and \ref{wnspsSEC} should suffice to understand the algebraic construction of the product of noises in Section \ref{wnSEC}.

In Section \ref{CPD} we repeat some results about units, CPD-semigroup and their generators and show that every spatial product system has a unique maximal completely spatial subsystem. In Section \ref{tpuSEC} we provide certain \it{geometric} operations and a \it{Trotter product} among units, that help in Section \ref{prodSEC} to understand the units in the product of spatial product systems and in Section \ref{toSEC} to show that completely spatial product systems are time ordered (allowing the definition of the index) and to show that the index is additive under our product. In Section \ref{wnSEC} we extend all results from spatial product systems to noises in a compatible way.

Section \ref{freeflow} is dedicated to a detailed analysis of \it{free flows}, that is time shifts on full Fock modules. We arrive at analogues to Fowler's [\refcite{Fow95}] results for the Hilbert space case, but it seems that our description is more specific, the derivation of the product system is more direct and the combinatorical problems of the time ordered Fock module are much simpler than those of the symmetric Fock space. In Section \ref{openSEC} we pose a couple of open and, we think, interesting problems. A good deal of them is motivated by the fact that the free flows in Section \ref{freeflow} can be described conveniently also by \it{free product systems}.

\lf\noindent
\bf{Conventions.~}
Throughout these notes we use a couple of results from [\refcite{BhSk00,LiSk01,BBLS04,Ske02}] which can also be found in [\refcite{Ske01}]. We reference, usually, to [\refcite{Ske01}] which is accessible through the author's homepage. [\refcite{Ske01}] contains also a detailed introduction to Hilbert modules. A less specific (for our purposes) reference for Hilbert modules is the book of Lance [\refcite{Lan95}]. A short introduction to Hilbert modules as we need them (only with few proofs) and everything about inductive limits we will be using can be found in [\refcite{BhSk00}].

Here we recall only that Hilbert \nbd{\cB}modules are always right $\cB$--mod\-ules with a \nbd{\cB}valued inner product. Here the \nbd{C^*}algebra $\cB$ will always be unital. A \hl{Hilbert \nbd{\cB}\nbd{\cB}module} (or \hl{two-sided} Hilbert module) is a Hilbert \nbd{\cB}module with a \hl{unital} representation of $\cB$ by adjointable (and, therefore, bounded and right linear) mappings.\footnote{In more recent articles, instead of \it{two-sided module}, we switched to the more standard name \it{correspondence}.} The algebra of all bounded adjointable mappings on a Hilbert module $E$ is denoted by $\sB^a(E)$. Whereas, $\sB^{a,bil}(E)$ denotes the subalgebra of \nbd{\cB}\nbd{\cB}linear (or two-sided) mappings. By $xy^*$ we denote the \hl{rank-one operator} $z\mapsto x\AB{y,z}$. The inner product on the \hl{tensor product} $E\odot F$ of two Hilbert \nbd{\cB}\nbd{\cB}modules $E,F$ is defined by $\AB{x\odot y,x'\odot y'}=\AB{y,\AB{x,x'}y'}$. Often we use only pre-Hilbert modules (a Hilbert module, except that it need not be complete). Constructions like tensor products $\odot,\otimes$ (where the latter is always that of vector spaces) and $\oplus$ are always understood algebraic, whereas we indicate completions by $\sbar{\odot}$, and so on.

\lf\noindent
\bf{A note on time lines.~}
These notes have been published first as Volterra Preprint 458, Centro Vito Volterra, University of Rome II, in March 2001, containing all definitions and results except Sections \ref{freeflow} and \ref{openSEC}. (Section \ref{freeflow} has been added in the end of 2001 and Section \ref{openSEC} in 2003.) As the inclusion of references to recent new results has caused more than once misunderstandings regarding priority of publications, in this revision we decided to put every reference to articles that have been written after these notes into footnotes. (Only Section \ref{openSEC} that, othewise would result in a single huge footnote, remains as it is.) These notes are the first place where spatial product systems of Hilbert modules have been defined, and where their product has been constructed. All other papers of which the author is (co-)author do refer to the present notes as primary source for spatial product systems and their product.

Also, as compared with earlier versions, proofs preceding their theorem have switched order with the theorem in order to underline that the paper contains proofs, and statements hidden in remarks became theorems and propositions in order to underline that the paper contains statements. Also some details left out in proofs for the reader and, thus, leading to quite condensed proofs, have now been filled in.

\bf{A note on terminology and title.~}
When first published we used the term \it{white noise} for what in Powers' terminology would be a spatial \it{\nbd{E_0}semigroup in standard form}. Following criticisms by L.\ Accardi and by C.\ K\"ostler we agree on that the terminology \it{white noise} is highly missleading. Actually, what we are considering corresponds very well to what Tsirelson introduced as \it{noise} replacing classical independence with \it{amalgamated monotone independence}; see Skeide [\refcite{Ske04}]. (We do not follow Hellmich, K\"ostler and K\"ummerer [\refcite{HKK04p}] who would rather say \it{Bernoulli shift}.) Throughout this revision \it{white noise} has been substituted with \it{noise}. An exception is the title where we put \it{white} into parenthesis in order to not change the title too much (obscuring the time lines and giving again rise to discussions about priority).
}


\section{Product systems, units and spatial product systems}\label{psSEC}

In this section we do not much more than recalling a few definitions and results from [\refcite{BhSk00,Ske02,BBLS04}] (that all may be found also in [\refcite{Ske01}]) and we define \it{spatial} product systems and their morphisms. This lays the basis for the product of spatial product systems in Section \ref{prodSEC}.

Let $\cB$ be a unital \nbd{C^*}algebra. A \hl{tensor product system of pre-Hilbert modules}, or for short a \hl{product system}, is a family $E^\odot=\bfam{E_t}_{t\in\R_+}$ of pre-Hilbert \nbd{\cB}\nbd{\cB}modules $E_t$ with a family of two-sided unitaries $u_{st}\colon E_s\odot E_t\rightarrow E_{s+t}$ $(s,t\in\R_+)$, fulfilling the associativity condition
\beq{\label{PSass}
\parbox{8cm}{
\xymatrix{
&E_r\odot E_s\odot E_t	\ar[dl]_{u_{rs}\odot\id_{E_t}}	\ar[dr]^{\id_{E_r}\odot u_{st}}	&
\\
E_{r+s}\odot E_t	\ar[dr]_{u_{(r+s)t}}	&	&E_r\odot E_{s+t}	\ar[dl]^{u_{r(s+t)}}
\\
&E_{r+s+t}&
}
}
}\eeq
where $E_0=\cB$ and $u_{s0}$ and $u_{0t}$ are the canonical identifications $x_s\odot b=x_sb$ and $b\odot x_t=bx_t$, respectively. Once the choice of $u_{st}$ is fixed, we always use the identification
\beq{ \label{tpsid}
E_s\odot E_t
~=~
E_{s+t}.
}\eeq
We speak of a \hl{tensor product system} of \hl{Hilbert modules} $E^{\sbar{\odot}}$, if $E_s\sbar{\odot} E_t=E_{s+t}$. We do not discuss the obvious generalizations to von Neumann modules.

A \hl{product subsystem} is a family ${E'}^\odot=\bfam{E'_t}_{t\in\R_+}$ of \nbd{\cB}\nbd{\cB}submodules $E'_t$ of $E_t$ such that $E'_s\odot E'_t=E'_{s+t}$ by restriction of the identification \eqref{tpsid}.

By the \hl{trivial} product system we mean $\bfam{\cB}_{t\in\R_+}$ where $\cB$ is equipped with its natural $\cB$--$\cB$--mod\-ule structure and inner product $\AB{b,b'}=b^*b'$.

A \hl{morphism} between product systems $E^\odot$ and $F^\odot$ is a family $w^\odot=\bfam{w_t}_{t\in\R_+}$ of mappings $w_t\in\sB^{a,bil}(E_t,F_t)$, fulfilling
\beq{ \label{psm}
w_{s+t}
~=~
w_s\odot w_t
}\eeq
and $w_0=\id_\cB$. A morphism is \hl{unitary}, \hl{contractive}, and so on, if $w_t$ is for all $t\in\R_+$. An \hl{isomorphism} of product systems is a unitary morphism. In rare occasions when we do not require the $w_t$ to be bounded (but still adjointable and two-sided) we speak of \hl{(possibly unbounded)} morphisms (necessarily of product systems of pre-Hilbert modules, because adjointable mappings between Hilbert modules are bounded, automatically).

Observe that, in general, there need not exist a projection morphism onto a subsystem.

A \hl{unit} for a product system $E^\odot=\bfam{E_t}_{t\in\R_+}$ is a family $\xi^\odot=\bfam{\xi_t}_{t\in\R_+}$ of elements $\xi_t\in E_t$ such that
\beq{ \label{unitp}
\xi_s\odot\xi_t
~=~
\xi_{s+t}
}\eeq
in the identification \eqref{tpsid} and $\xi_0=\U\in\cB=E_0$. By $\sU(E^\odot)$ we denote the set of all units for $E^\odot$. A unit $\xi^\odot$ is \hl{unital} and \hl{contractive}, if $\AB{\xi_t,\xi_t}=\U$ and $\AB{\xi_t,\xi_t}\le\U$, respectively. A unit is \hl{central}, if $b\xi_t=\xi_tb$ for all $t\in\R_+,b\in\cB$.

Obviously, a morphism $w^\odot\colon E^\odot\rightarrow F^\odot$ sends units to units. For this the requirement $w_0=\id_\cB$ is necessary. For a subset $S\subset\sU(E^\odot)$ of units for $E^\odot$ we denote by $w^\odot S\subset\sU(F^\odot)$ the subset of units for $F^\odot$, consisting of the units $w\xi^\odot=\bfam{w_t\xi_t}_{t\in\R_+}$ $(\xi^\odot\in S)$.

Now we are ready to define spatial product systems.

\bdefi \label{spatdef}
A \hl{spatial} product system is a pair $(E^\odot,\om^\odot)$ consisting of a product system $E^\odot$ and a central unital unit $\om^\odot$, the \hl{reference unit}. A \hl{spatial subsystem} of a spatial product system $(E^\odot,\om^\odot)$ is a subsystem that contains $\om^\odot$ and that is spatial with $\om^\odot$ as reference unit. A morphism $w^\odot$ between spatial product systems is \hl{spatial}, if both $w^\odot$ and ${w^*}^\odot$ preserve the reference units.
\edefi

More loosely, we speak of a spatial product system $E^\odot$, if we can turn it into a spatial one by choosing a central unital unit $\om^\odot$. But, we must be aware that structures derived from that unit (as, for instance, the product of spatial systems in Section \ref{prodSEC}) will depend on the choice of $\om^\odot$.\footnote{The reference unit can also be used to pose measurability or continuity conditions on the product system and one can show that these do not depend on the choice of the reference unit (as long as the two reference units are continuous among themselves in the sense of Lemma \ref{contlem}; see [\refcite{Ske03b}]). In the case of Arveson systems, we obtain back Arveson's measurability conditions.} However, even for Arveson systems\footnote{which are known to be isomorphic if they are algebraically isomorphic; see Liebscher [\refcite{Lie00p1}]} it was a long time unclear, in how far spatial Arveson systems are isomorphic spatial product systems in our sense, if they are different only for the choice of a different reference unit.\footnote{The related problem, open for a long time, is the question whether the automorphism group of an Arveson system acts transitively on its set of units. Only recently Tsirelson [\refcite{Tsi04p1}] provided us with a concrete counter example where the automorphism group of a type II$_1$ Arvesons system does not act transitively on the set of normalized units. So the spatial structure of a product system, indeed, may depend on the choice of the reference unit.}

\lf
We close this section with some remarks.

It is known that spatial Arveson systems are those that have units, and because every unit of an Arveson system is central the definitions coincide. On the other hand, it is easy to write down examples of product systems that have units but none of these is central; see [\refcite{BBLS04}] (or [\refcite{Ske03b}]). So, why do we call spatial those product systems that admit a central unital unit? The answer lies in the close relationship between Arveson systems and \nbd{E_0}semigroups (semigroups of unital endomorphisms) on $\sB(H)$ [\refcite{Arv89,Arv89a}] and Powers' original definition of spatial \nbd{E_0}semigroups in [\refcite{Pow88}]. To every \nbd{E_0}semigroup there is an associated product system (we discuss the Hilbert module version in Section \ref{wnspsSEC}) and it is easy to show that the so-called \it{intertwining semigroups of isometries} for the \nbd{E_0}semigroup correspond (one-to-one if $E$ is full) to central unital units of the associated product system. (Notice that this fails, if $\cB$ is nonunital. This is one of the main reasons, why we stick to unital $\cB$.) And Powers' definition says an \nbd{E_0}semigroup is spatial, if it admits an intertwining semigroup of isometries.

There are plenty of spatial product systems. In Section \ref{toSEC} we will see that the subclass of \it{completely spatial} product systems of Hilbert modules consists precisely of the time ordered product systems and by Corollary \ref{maxspat} below any spatial product system contains a (unique) maximal completely spatial subsystem.\footnote{A construction by Liebscher [\refcite{Lie00p1}] allows to construct from every Arveson system loads of spatial Arveson systems that are not completely spatial. We believe that this can be done for arbitrary (continuous in the sense of [\refcite{Ske03b}]) spatial product systems of Hilbert modules, thus, providing us with many examples of highly nontrivial spatial product systems.} By a result from [\refcite{BBLS04}] ([\refcite{Ske01}, Corollary 13.2.13]) type $\text{I}$ product systems of von Neumann modules are (strong closures of) time ordered systems. Therefore, if a product system of von Neumann modules contains a single continuous unit $\xi^\odot$, then the subsystem generated by $\xi^\odot$ is time ordered and, therefore, contains a central unital unit, namely, its vacuum unit (see Section \ref{toSEC}). In other words, in the context of von Neumann modules existence of a single continuous unit is sufficient to know that a product system is spatial. One can show (see [\refcite{BBLS04}]) that this is equivalent to the result by Christensen and Evans [\refcite{ChrEv79}] that (rephrased suitably) bounded derivations with values in a von Neumann module are inner.

\section{Units and CPD-semigroups}\label{CPD}

The notions introduced so far, are sufficient to understand the constructions of the products in Sections \ref{prodSEC} and \ref{wnSEC}, and it is possible to read them now. These constructions are completely algebraic and extend by well-known compatibility conditions to any desired completion (or closure in the case of von Neumann modules). In this section we recall the basic classification of product systems by units (mainly from [\refcite{BBLS04}]) and we draw first consequences from existence of a central unital reference unit. In particular, we define completely spatial product systems and show that every spatial product system contains a unique maximal completely spatial subsystem.

A crucial role in the analysis of type I Arveson systems is played by a semigroup of positive definite kernels on the set of units defined by $(u^\otimes,{u'}^\otimes)\mapsto\AB{u_t,u'_t}$. In [\refcite{Arv89}], the generator of this semigroup is named the \it{covariance function} of an Arveson system. Also for product systems of Hilbert modules the inner products of units determine a semigroup of kernels, however, the more noncommutative structure of \nbd{\cB}\nbd{\cB}modules (even, or actually, in particular, if $\cB$ is commutative) where $bx=xb$ happens only rarely, forces us to consider the mappings $b\mapsto\AB{\xi_t,b\xi'_t}$ rather than the matrix elements $\AB{\xi_t,\xi'_t}$.

For us a \hl{kernel} on a set $S$ is a mapping $\eK\colon S\times S\rightarrow\sB(\cB)$ into the bounded mappings on $\cB$. (This, clearly, contains the well-known notion of \nbd{\C}valued kernels, if we consider an element $w\in\C$ as mapping $z\mapsto wz$ in $\sB(\C)$.) According to the definition in [\refcite{BBLS04}] a kernel $\eK$ is \hl{completely positive definite}, if
\beqn{
\sum_{i,j}b_i^*\eK^{\sigma_i,\sigma_j}(a_i^*a_j)b_j
~\ge~
0
}\eeqn
for all choices of finitely many $\sigma_i\in S;a_i,b_i\in\cB$. (Our definition in [\refcite{BBLS04}] was inspired very much by that in Accardi and Kozyrev [\refcite{AcKo01}]. We emphasize, however, that the definition in [\refcite{AcKo01}] is weaker, but due to additional structure present in their concrete problem also their kernel is completely positive definte in our sense; see [\refcite{Ske01}, Lemmata 5.2.7 and 5.3.5].) Clearly, any kernel $\eK$ of the form $\eK^{\sigma,\sigma'}(b)=\AB{x_\sigma,bx_{\sigma'}}$ for some elements $x_\sigma$ $(\sigma\in S)$ in some pre--Hilbert \nbd{\cB}\nbd{\cB}module $E$ is completely positive definite. Moreover, any completely positive definite kernel can be recovered in that way by its \hl{Kolmogorov decomposition}; see [\refcite{BBLS04}] and [\refcite{Ske01}, Theorem 5.2.3].

The family $\eU=\bfam{\eU_t}_{t\in\R_+}$ of kernels $\eU_t$ on $\sU(E^\odot)$, defined by setting
\beq{ \label{unisg}
\eU_t^{\xi,\xi'}(b)
~=~
\AB{\xi_t,b\xi'_t}
}\eeq
is a semigroup under pointwise composition of kernels, as
\beqn{
\eU_{s+t}^{\xi,\xi'}(b)
~=~
\AB{\xi_{s+t},b\xi'_{s+t}}
~=~
\AB{\xi_s\odot\xi_t,b\xi'_s\odot\xi'_t}
~=~
\bAB{\xi_t,\AB{\xi_s,b\xi'_s}\xi'_t}
~=~
\eU_t^{\xi,\xi'}\circ\eU_s^{\xi,\xi'}(b),
}\eeqn
and all $\eU_t$ are comletely positive definite. We say $\eU$ is the \hl{CPD-semigroup} associated with the product system $E^\odot$.

Every CPD-semigroup, i.e\, in particular, every CP-semigroup, can be recovered in this way from its \hl{GNS-system}; see [\refcite{BBLS04}] and [\refcite{Ske01}, Theorem 11.3.5]. In other words, any CPD-semigroup is obtained from units of a product system as in \eqref{unisg}. However, the converse need not be true as follows from the existence of nonspatial Arveson systems. Nevertheless, by [\refcite{Ske01}, Proposition 11.2.4] any subset $S\subset\sU(E^\odot)$ of units of a product system $E^\odot$ \hl{generates} a product subsystem ${E^S}^\odot=\bfam{E^S_t}_{t\in\R_+}$ consisting of the spaces
\beq{\label{tensgen}
E^S_t
~=~
\ls\bset[b_n\xi^n_{t_n}\odot\ldots\odot b_1\xi^1_{t_1}b_0]{n\in\N,b_i\in\cB,{\xi^i}^\odot\in S,t_n+\ldots+t_1=t}.
}\eeq
(Compare the definition of the lattice $\bJ_t$ in the beginning of Section \ref{prodSEC}.)

The CPD-semigroup $\eU\upharpoonright S$ is \hl{uniformly continuous}, if the semigroups in \eqref{unisg} are uniformly continuous for all $\xi^\odot,{\xi'}^\odot\in S$. In this case we say $S$ is a \hl{continuous} subset of units. In particular, a single unit is \hl{continuous}, if the subset $\CB{\xi^\odot}$ is continuous.

\bob
Obviously, every central unital unit is continuous. And if a central unit $\om^\odot$ is continuous, then we may modify it to be unital by normalizing it to $\om_t\sqrt{\AB{\om_t,\om_t}}^{-1}$, because $\AB{\om_t,\om_t}$ is a continuous semigroup of positive central elements in $\cB$.
\eob

The type of a product sytem was defined in [\refcite{BBLS04}], in analogy with that of an Arveson system, indicating in how far the product system is generated by its units. The only difference is that Hilbert modules have several different topologies and topology enters in two essentially different ways. Firstly, there are different closures in which the product system might be generated by its units and, secondly, there are different topologies in which a CPD-semigroup might be continuous. We repeat here only the relevant part of the definitions from [\refcite{BBLS04}], that is we consider only sets of units that lead to uniformly continuous CPD-semigroups and we do not give the version for von Neumann modules.

A product system $E^\odot=\bfam{E_t}_{t\in\R_+}$ of pre-Hilbert modules is of \hl{type $\ul{\text{I}}$}, if it is \hl{generated} by some continuous set $S\subset\sU(E^\odot)$ of units, i.e.\ if $E^\odot={E^S}^\odot$. It is of \hl{type I}, if $E^\odot$ (or $E^{\sbar{\odot}}$ in the case of Hilbert modules) is the norm closure of ${E^S}^\odot$. We say the set $S$ is \hl{generating}.

Now we are ready to define completely spatial product systems.

\bdefi
A spatial product system $(E^\odot,\om^\odot)$ is \hl{completely spatial} of type $\ul{\text{I}}$, and so on, if it is type $\ul{\text{I}}$, and so on, such that the subset $S\subset\sU(E^\odot)$ of units making it type $\ul{\text{I}}$, and so on, can be chosen such that $\om^\odot\in S$. If we do not specify the type of completely spatial, then we mean always type $\text{I}$.
\edefi

The following is a (slightly weaker) reformulation of [\refcite{BBLS04}, Theorem 4.4.12] (or [\refcite{Ske01}, Lemma 11.6.6]).\footnote{It has a much stronger counterpart in the discussion of \it{continuous} product systems in [\refcite{Ske03b}] where $\om^\odot$ may be an arbitrary continuous unit and $S$ a set of general \it{continuous sections} of $E^{\sbar{\odot}}$.}

\blem \label{contlem}
Let $S$ be a subset of units in a spatial product system $(E^{\sbar{\odot}},\om^\odot)$ of Hilbert modules. Then $\eU\upharpoonright\CB{\om^\odot}\cup S$ is uniformly continuous, if and only $\eU\upharpoonright\CB{\om^\odot,\xi^\odot}$ is uniformly continuous for all $\xi^\odot\in S$.
\elem

\bcor \label{maxspat}
The set
\beqn{
\sU_\om(E^{\sbar{\odot}})
~:=~
\bCB{\xi^\odot\in\sU(E^{\sbar{\odot}})\colon\eU\upharpoonright\CB{\om^\odot,\xi^\odot}\text{\rm{~is uniformly continuous}}}
}\eeqn
is the maximal subset $S$ of $\sU(E^{\sbar{\odot}})$ containing $\om^\odot$ for which $\eU\upharpoonright S$ is uniformly continuous. Consequently, each spatial product system of Hilbert modules has a unique \hl{maximal completely spatial subsystem} $({E^{\sU_\om}}^{\sbar{\odot}},\om^\odot)$.
\ecor

\proof
By Lemma \ref{contlem} $S$ is continuous. Moreover, every other continuous subset $S'$ containing $\om^\odot$, clearly, fulfills the condition in Lemma \ref{contlem} and, therefore, by defintion is contained in $S$.~\qed

\section{The Trotter product of units}\label{tpuSEC}

In this section we provide some ``geometric'' operations among units. First, we construct a sort of \it{artithmetic mean} for (continuous) units in arbitrary product systems (although here we will prove only the simpler spatial case). Then, applying results from [\refcite{BBLS04}] to the case of spatial product systems, we use the mean to construct a \it{Trotter product} of units. The \it{Trotter product} of units will show us how to \it{compose} units from different factors in the product of spatial product systems discussed in Section \ref{prodSEC}. It will help us to show that the index is \it{additive}.

Many properties of the units and the operations among them have more concrete interpretations, when applied to units in a time-ordered product system as discussed in Section \ref{toSEC}. For instance, what we define to be an \it{exponential unit} in a spatial product system, really, corresponds to an exponential unit in a time-ordered product system. However, in particular the arithmetic mean operation, which is defined for units in not necessarily spatial product systems, would make it necessary to refer to the deep embeddability result (into a time-ordered product system of von Neumann modules only) in [\refcite{BBLS04}]. As we wish to avoid the discussion of von Neumann modules, we have keep the discussion of these properties on an abstract level.

\lf
Let $E^{\sbar{\odot}}$ be a product system of Hilbert modules with a continuous subset of units $S$ and dentote by $\eL_S=\frac{d}{dt}\eU\big|_{t=0}\upharpoonright S$ the generator of the CPD-semigroup $\eU\upharpoonright S$. 

\blem \label{Tlem}
Let ${\xi^\ell}^\odot\in S$ $(\ell=1,2)$. Then for all $\vk^1,\vk^2\in\C$ with $\vk^1+\vk^2=1$ the limit
\beq{\label{Trotter}
\xi_t
~=~
\lim_{n\to\infty}\bfam{\vk^1\xi^1_{\frac{t}{n}}+\vk^2\xi^2_{\frac{t}{n}}}^{\odot n}
}\eeq
exists in norm, $\xi^\odot=\bfam{\xi_t}_{t\in\R_+}$ is a unit, too, and the set $S\cup\CB{\xi^\odot}$ is still continuous. Moreover, for all ${\xi'}^\odot\in S\cup\CB{\xi^\odot}$ we have $\eL^{\xi',\xi}=\vk^1\eL^{\xi',\xi^1}+\vk^2\eL^{\xi',\xi^2}$.
\elem

We will prove Lemma \ref{Tlem} and a generalization to nets, which we need in Section \ref{prodSEC}, in the appendix and only for spatial product systems. A full proof would require to repeat a good deal more from [\refcite{BhSk00,BBLS04}] will appear elswhere. (See  Liebscher and Skeide [\refcite{LiSk05p}].)

Observe that there is at most one unit $\xi^\odot$ fulfilling $\eL^{\xi',\xi}=\vk^1\eL^{\xi',\xi^1}+\vk^2\eL^{\xi',\xi^2}$ for all ${\xi'}^\odot\in S$, because $\xi^\odot$ is contained in the type I subsystem generated by $\CB{{\xi^1}^\odot,{\xi^2}^\odot}$ and the inner products within this subsystem are determined completely by the generator $\eL$.

The lemma has an obvious generalization to $n$ units ${\xi^1}^\odot,\ldots,{\xi^n}^\odot$ with $n$ complex numbers $\vk^1+\ldots+\vk^n=1$. We use the notation $\bfam{\vk^1\xi^1\bplus\ldots\bplus\vk^n\xi^n}^\odot$. Like the arithmetic mean, the operation $\bplus$ is commutative. It is associative in the sense that
\bmun{
\Bfam{(\vk+\vk')\Bfam{\frac{\vk}{\vk+\vk'}\xi\bplus\frac{\vk'}{\vk+\vk'}\xi'}\bplus\vk''\xi''}^\odot
\\
~=~
\bfam{\vk\xi\bplus\vk'\xi'\bplus\vk''\xi''}^\odot
\\
~=~
\Bfam{\vk\xi\bplus(\vk'+\vk'')\Bfam{\frac{\vk'}{\vk'+\vk''}\xi'\bplus\frac{\vk''}{\vk'+\vk''}\xi''}}^\odot
}\emun
$(\vk+\vk'\ne0\ne\vk'+\vk'')$. To see this, just look at the the generator and use uniqueness.

\lf
Now we pass to a spatial product system $(E^{\sbar{\odot}},\om^\odot)$ and put $S=\sU_\om(E^{\sbar{\odot}})$. So $\eL$ is now the generator of the CPD-semigroup $\eU\upharpoonright\sU_\om(E^{\sbar{\odot}})$. Since $\eU^{\om,\xi}_t(b)=b\eU^{\om,\xi}_t(\U)$, the elements $\eU^{\om,\xi}_t(\U)$ form a semigroup in $\cB$. For $\xi^\odot\in\sU_\om(E^{\sbar{\odot}})$ we denote by $\beta_\xi=\eL^{\om,\xi}(\U)$ the generator of this semigroup. The crucial result [\refcite{BBLS04}, Theorem 5.1.2] (or [\refcite{Ske01}, Theorem 13.1.2]) asserts that the kernel $\eL_0$ defined by
\beq{\label{msgen}
\eL_0^{\xi,\xi'}(b)
~=~
\eL^{\xi,\xi'}(b)-\beta_\xi^*b-b\beta_{\xi'}
}\eeq
is completely positive definite. We say a kernel $\eL$ allowing for a decomposition $\eL^{\xi,\xi'}(b)=\eL_0^{\xi,\xi'}(b)+\beta_\xi^*b+b\beta_{\xi'}$ for a completely positive definite kernel $\eL_0$ and suitable $\beta_\xi\in\cB$ has \hl{CE-form} (Chri\-sten\-sen-Evans form) or is a \hl{CE-generator}.

Observe that $\beta_\om=0$ and, therefore, $\eL_0^{\om,\xi}=0=\eL_0^{\xi,\om}$ for all $\xi^\odot$. More generally, if for $\beta\in\cB$ we denote by ${\om^\beta}^\odot$ the unit $\bfam{\om_te^{t\beta}}_{t\in\R_+}$, then $\beta_{\om^\beta}=\beta$ and $\eL_0^{\om^\beta,\xi}=0=\eL_0^{\xi,\om^\beta}$. Note that these two properties determine ${\om^\beta}^\odot$ uniquely.

\bprop\label{utoe}
Let $\xi^\odot\in\sU_\om(E^{\sbar{\odot}})$ and $\beta\in\cB$. Then ${\zeta}^\odot=\bfam{\xi\bplus\om^{\beta}\bplus-\om}^\odot$ is the unique unit in $\sU_\om(E^{\sbar{\odot}})$ satisfying $\beta_{\zeta}=\beta_\xi+\beta$ and $\eL_0^{\xi',\zeta}=\eL_0^{\xi',\xi}$ for all $\xi'^\odot\in\sU_\om(E^{\sbar{\odot}})$.
\eprop

\proof
Uniqueness is clear as the stated properties determin the generator. By the three term version of Lemma \ref{Tlem} and the preceding computations we find
\beqn{
\beta_\zeta
~=~
\eL^{\om,\zeta}(\U)
~=~
\eL^{\om,\xi}(\U)+\eL^{\om,\om^{\beta}}(\U)-\eL^{\om,\om}(\U)
~=~
\beta_\xi+\beta-0
}\eeqn
and
\bmun{
\eL_0^{\xi',\zeta}(b)
\\
~=~
\eL^{\xi',\zeta}(b)-\beta_{\xi'}^*b-b\beta_\zeta
~=~
\eL^{\xi',\xi}(b)+\eL^{\xi',\om^{\beta}}(b)-\eL^{\xi',\om}(b)-\beta_{\xi'}^*b-b(\beta_\xi+\beta)
\\
~=~
\eL^{\xi',\xi}(b)+(0+\beta_{\xi'}^*b+b\beta)-(0+\beta_{\xi'}^*b+b0)-\beta_{\xi'}^*b-b(\beta_\xi+\beta)
\\
~=~
\eL_0^{\xi',\xi}(b).
}\emun
\qed

\lf
Proposition \ref{utoe} tells us that we may always pass from a continuous unit to a unit with $\beta_\xi=0$ without changing the CPD-part $\eL_0$ of the generator. We call a unit with $\beta_\xi=0$ an \hl{exponential unit}, because in Section \ref{toSEC} such units will show to be precisely those that consist of exponential vectors.

\bcor\label{etou}
Let $\xi^\odot$ be a continuous unit. Then $\zeta^\odot=\bfam{\xi\bplus\om^{-\beta_\xi}\bplus-\om}^\odot$ by Proposition \ref{utoe} is the unique exponential unit such that $\xi^\odot=\bfam{\zeta\bplus\om^{\beta_\xi}\bplus-\om}^\odot$.
\ecor

\bdefi
Let ${\xi^\ell}^\odot$ $(\ell=1,2)$ be continuous units in a spatial product system $(E^{\sbar{\odot}},\om^\odot)$. By the \hl{Trotter product} of ${\xi^1}^\odot$ and ${\xi^2}^\odot$ we mean the unit
\beqn{
\bfam{\xi^1\circledcirc\xi^2}^\odot
~:=~
\bfam{\xi^1\bplus\xi^2\bplus-\om}^\odot
}\eeqn
\edefi

In this notation Corollary \ref{etou} reads $\xi^\odot=\bfam{\zeta\circledcirc\om^{\beta_\xi}}^\odot$ where $\zeta^\odot=\bfam{\xi\circledcirc\om^{-\beta_\xi}}^\odot$ is the unique exponential unit not changing the CPD-part of the generator when substituting $\xi^\odot$ with $\zeta^\odot$.

\bprop
The Trotter product is associative.
\eprop

\proof
As may be checked by looking at the generators, both bracketings $\bfam{(\xi^1\circledcirc\xi^2)\circledcirc\xi^3}^\odot$ and $\bfam{\xi^1\circledcirc(\xi^2\circledcirc\xi^3)}^\odot$ lead to the same expression $\bfam{\xi^1\bplus\xi^2\bplus\xi^3\bplus-2\om}^\odot$.~\qed

\bcor
Let ${\xi^\ell}^\odot$ $(\ell=1,2)$ be continuous units and denote by ${\zeta^\ell}^\odot$ the corresponding exponential units according to Corollary \ref{utoe}. Then $\bfam{\xi^1\circledcirc\xi^2}^\odot=\bfam{(\zeta^1\circledcirc\zeta^2)\circledcirc(\om^{\beta_{\xi^1}}\circledcirc\om^{\beta_{\xi^2}})}^\odot$, where $\bfam{\zeta^1\circledcirc\zeta^2}^\odot$ is the unique exponential unit fulfilling  $\eL_0^{\xi',\zeta^1\circledcirc\zeta^2}=\eL_0^{\xi',\xi^1\circledcirc\xi^2}=\eL_0^{\xi',\xi^1}+\eL_0^{\xi',\xi^2}$ and where $\bfam{\om^{\beta_{\xi^1}}\circledcirc\om^{\beta_{\xi^2}}}^\odot=\om^{\beta_{\xi^1}+\beta_{\xi^2}}$.
\ecor

A warning: Observe that $\eL^{\xi',\xi^1\circledcirc\xi^2}=\eL^{\xi',\xi^1}+\eL^{\xi',\xi^2}-\eL^{\xi',\om}\ne\eL^{\xi',\xi^1}+\eL^{\xi',\xi^2}$.

\section{The product of spatial product systems} \label{prodSEC}

In this section we construct the \it{product} of spatial product systems. Unlike the tensor product of Arveson systems, our product works for Hilbert modules over arbitrary \nbd{C^*}algebras. When applied to Arveson systems it coincides only occasionally with the tensor product, but in general it is only a subsystem of the tensor product. Using existence of projection morphisms onto the factors, we establish a universal property of our product (in general not shared by the tensor product of Arveson systems). We give a  complete characterization of the continuous units of the product in terms of the continuous units of the factors.

Set $\bJ_t=\bCB{\et=(t_n,\ldots,t_1)\in(0,t]^n\colon n\in\N,\abs{\et}=t_n+\ldots+t_1=t}$. We consider $\bJ_t$ as a lattice taking its order structure from the set of interval partitions $\bI_t=\CB{\es=(s_n,\ldots,s_1)\in(0,t]^n\colon n\in\N,t=s_n>\ldots>s_1>s_0=0}$ via the bijection $\es\mapsto(s_n-s_{n-1},\ldots,s_1-s_0)$.

Let ${E^\ell}^\odot$ $(\ell=1,2)$ be two product systems of pre-Hilbert modules. The idea to construct a product of ${E^1}^\odot$ and ${E^2}^\odot$ is to think of a space spanned by elements of the form
\beq{ \label{genel}
x_{t_n}^n\odot\ldots\odot x_{t_1}^1
}\eeq
for $\et\in\bJ_t$ and $x_{t_i}^i$ either in $E_{t_i}^1$ or in $E_{t_i}^2$. There is no problem to provide such a space as vector space by an inductive limit over $\bJ_t$. (Every $x_{t_i}^i\in E^{\ell_i}_{t_i}$ can be written as an element in the tensor product of several $E^{\ell_i}_s$ to smaller times $s$, thus, giving rise to a refinement of $\et$.) Since we want to have a two-sided pre-Hilbert module, we have to face the problem to define an inner product of elements $x^1\in E_{t_i}^1$ and $x^2\in E_{t_i}^2$.

A first attempt could be to think about the direct sum, i.e.\ $\AB{x^1,x^2}=0$. This is, indeed, possible and results into a product system. Even the units of each ${E^\ell}^\odot$ embed as units into the new product system. However, $\eU^{\xi^1,\xi^2}_t$ is $0$ for $t>0$ and, of course, $\id_\cB$ for $t=0$ so that the semigoup $\eU^{\xi^1,\xi^2}$ will not be continuous in any reasonable topology. Therefore, even if the factors are type $\text{I}$, a product constructed in that way will never be type I. It is also not difficult to see that already in the Hilbert space case the new product system, in general, does not consist of separable Hilbert spaces.

A solution of this difficulty is possible, if we restrict to spatial product systems $({E^\ell}^\odot,{\om^\ell}^\odot)$ and require that the images of both reference units are no longer distinguished in the product and serve there as a reference unit $\om^\odot$. In other words, we have to identify at least the vectors $\om^1_t$ and $\om^2_t$ in the new product systems, while only the complements of $\om^1_t$ and $\om^2_t$ remain orthogonal. That is, we require that inner products of elements $x^1\in E^1_t,x^2\in E^2_t$ from different factors have the form
\beq{\label{ipdef}
\AB{x^1,x^2}
~=~
\AB{x^1,\om^1_t}\AB{\om^2_t,x^2}.
}\eeq
The follwoing theorem is proved by making these ideas precise in an inductive limit construction.

\bthm\label{prodpropthm}
There exists a spatial product system $(F^\odot,\om^\odot)$ fulfilling the following properties.

\begin{enumerate}
\item
$(F^\odot,\om^\odot)$ contains $({E^\ell}^\odot,{\om^\ell}^\odot)$ $(\ell=1,2)$ as spatial subsystems and is generated by them, that is, $F_t$ is spanned by epressions like \eqref{genel}.

\item
The inner product of members $x^1\in E^1_t$ and $x^2\in E^2_t$ is given by \eqref{ipdef}.
\end{enumerate}
Moreover, every spatial product system fulfilling these properties is canonically isomorphic to $(F^\odot,\om^\odot)$.
\ethm

\bdefi\label{spatproddef}
We call $(F^\odot,\om^\odot)$ the \hl{product} of $({E^\ell}^\odot,{\om^\ell}^\odot)$ and we will denote it by $(\bfam{{E^1}\circledcirc{E^2}}^\odot,\om^\odot)$. By $(\bfam{{E^1}\circledcirc{E^2}}^{\sbar{\odot}},\om^\odot)$ we denote its completion which is a spatial product system of Hilbert modules.
\edefi

\proof[Proof of Theorem \ref{prodpropthm}.~]
Set $\wh{E}^\ell_t=(\U-\om^\ell_t{\om^\ell_t}^*)E^\ell_t$, so that $E^\ell_t=\om^\ell_t\cB\oplus\wh{E}^\ell_t$, and define
\beqn{
\breve{F}_t
~=~
\om_t\cB\oplus\wh{E}^1_t\oplus\wh{E}^2_t
}\eeqn
(where $\om_t\cB\cong\cB$ via $\om_t\mapsto\U$ is the one-dimensional two-sided \nbd{\cB}module). The inner product of elements $x^\ell\in E^\ell_t$ is now defined by identifying them via $E^\ell_t=\om^\ell_t\cB\oplus\wh{E}^\ell_t\cong\om_t\cB\oplus\wh{E}^\ell_t\subset \breve{F}_t$ as elements of $\breve{F}_t$. (In other words, $\AB{x^1,x^2}=\AB{x^1,\om^1_t}\AB{\om^2_t,x^2}$ for elements from different factors, while the inner product of elements from the same factor remains unchanged.)

Let us fix $s,t\in\R_+$. We have
\beqn{
E^\ell_{s+t}
~=~
(\om^\ell_s\cB\oplus\wh{E}^\ell_s)\odot(\om^\ell_t\cB\oplus\wh{E}^\ell_t)
~=~
\om^\ell_{s+t}\cB\,\oplus\,(\om^\ell_s\odot\wh{E}^\ell_t)\,\oplus\,(\wh{E}^\ell_s\odot\om^\ell_t)\,\oplus\,(\wh{E}^\ell_s\odot\wh{E}^\ell_t)
}\eeqn
so that $\wh{E}^\ell_{s+t}=(\om^\ell_s\odot\wh{E}^\ell_t)\,\oplus\,(\wh{E}^\ell_s\odot\om^\ell_t)\,\oplus\,(\wh{E}^\ell_s\odot\wh{E}^\ell_t)$, whence
\bmu{ \label{Fst}
\breve{F}_{s+t}
\\
~=~
\om_{s+t}\cB\,\oplus\,(\om^1_s\odot\wh{E}^1_t)\,\oplus\,(\wh{E}^1_s\odot\om^1_t)\,\oplus\,(\wh{E}^1_s\odot\wh{E}^1_t)\,\oplus\,(\om^2_s\odot\wh{E}^2_t)\,\oplus\,(\wh{E}^2_s\odot\om^2_t)\,\oplus\,(\wh{E}^2_s\odot\wh{E}^2_t).
\\
~
}\emu
On the other hand,
\bmun{
\breve{F}_s\odot\breve{F}_t
~=~
(\om_s\cB\,\oplus\,\wh{E}^1_s\,\oplus\,\wh{E}^2_s)\odot(\om_t\cB\,\oplus\,\wh{E}^1_t\,\oplus\,\wh{E}^2_t)
\\
~=~
\om_s\odot\om_t\cB\,\oplus\,(\om_s\odot\wh{E}^1_t)\,\oplus\,(\wh{E}^1_s\odot\om_t)\,\oplus\,(\wh{E}^1_s\odot\wh{E}^1_t)\,\oplus\,(\om_s\odot\wh{E}^2_t)\,\oplus\,(\wh{E}^2_s\odot\om_t)\,\oplus\,(\wh{E}^2_s\odot\wh{E}^2_t)
\\\,\oplus\,(\wh{E}^1_s\odot\wh{E}^2_t)\,\oplus\,(\wh{E}^2_s\odot\wh{E}^1_t). 
}\emun
In other words, sending in \eqref{Fst} $\om_{s+t}$ to $\om_s\odot\om_t$ and $\om^\ell_s$ and $\om^\ell_t$ to $\om_s$ and $\om_t$, respectively, defines an embedding of $\breve{F}_{s+t}$ into $\breve{F}_s\odot\breve{F}_t$ as a two-sided submodule. Clearly, for any $\et\in\bJ_t$ we may define a two-sided isometric embedding $\breve{F}_t\rightarrow\breve{F}_\et:=\breve{F}_{t_n}\odot\ldots\odot\breve{F}_{t_1}$ in a similar way. Finally, let $\et\ge\es\in\bJ_t$, i.e.\ $\et=(s^m_{\ell_m},\ldots,s^m_1,~\ldots~,s^1_{\ell_1},\ldots,s^1_1)$ (denoted as $\es_m\smallsmile\ldots\smallsmile\es_1$ in [\refcite{BhSk00,Ske01}]) where $\es_j=(s^j_{\ell_j},\ldots,s^j_1)\in\bJ_{s_j}$. By taking the tensor product of the mappings $\breve{F}_{s_j}\rightarrow\breve{F}_{\es_j}$ ($j=m,\ldots,1$) we end up with embeddings
\beqn{
\beta_{\et\es}\colon
\breve{F}_\es
~\rightarrow~
\breve{F}_\et.
}\eeqn
Such embeddings have been considered very carefully in [\refcite{BhSk00}] in the construction of product systems from CP-semigroups and in [\refcite{BBLS04}] in the construction of product systems from CPD-semigroups so that here we can proceed quickly (see [\refcite{Ske01}, Section 11.3]). See [\refcite{BhSk00}] or [\refcite{Ske01}, Appendix A.10] for inductive limits.

Clearly, $\beta_{\et\er}\beta_{\er\es}=\beta_{\et\es}$ for all $\et\ge\er\ge\es$ so that the $\breve{F}_\et$ with the $\beta_{\et\es}$ form an inductive system with an inductive limit $F_t$ which is a pre-Hilbert \nbd{\cB}\nbd{\cB}module. Indentifying $\breve{F}_\et$ with its image under the canonical embedding into $F_t$ (as prectised always in [\refcite{BhSk00,BBLS04,Ske01}]), we see that $F_t$, indeed, is spanned by elements of the form \eqref{genel} and the inner product of elements $x^i_{t_i},y^i_{t_i}$ is that of elements in $\breve{F}_{t_i}$. Setting $F_0=\cB$, as in the proof of [\refcite{BhSk00}, Theorem 4.8] we show that the family $F^\odot=\bfam{F_t}_{t\in\R_+}$ is a product system under the identification
\beq{\label{Ftpdef}
(x_{s_m}^m\odot\ldots\odot x_{s_1}^1)\odot(y_{t_n}^n\odot\ldots\odot y_{t_1}^1)
~=~
x_{s_m}^m\odot\ldots\odot x_{s_1}^1\odot y_{t_n}^n\odot\ldots\odot y_{t_1}^1
}\eeq
and $\om^\odot=\bfam{\om_t}_{t\in\R_+}$ is a unital central unit for $F^\odot$.

Uniqueness is obvious.~\qed

\bcor
Every unit ${\xi^\ell}^\odot$ for ${E^\ell}^\odot$ gives rise to a unit for $\bfam{{E^1}\circledcirc{E^2}}^\odot$ also denoted by ${\xi^\ell}^\odot$. In particular, we have ${\om^\ell}^\odot=\om^\odot$. The CPD-semigroup $\eU\upharpoonright\sU_\om({E^1}^\odot)\cup\sU_\om({E^2}^\odot)$ is uniformly continuous.
\ecor

\proof
The the embedded elements of a unit form again a unit follows from \eqref{Ftpdef}. Continuity of the units follows by Lemma \ref{contlem} because the new reference unit restricted to the subsystem is just the old reference unit.~\qed

\lf
Now we wish to specify $\sU_\om(\bfam{{E^1}\circledcirc{E^2}}^{\sbar{\odot}})$ better. For that goal we construct projection morphisms ${p^\ell}^\odot$ onto the subsystems ${E^\ell}^\odot$. The projection morphisms will, then, decompose a given unit $\xi^\odot$ for $\bfam{{E^1}\circledcirc{E^2}}^{\sbar{\odot}}$ into components in ${E^\ell}^{\sbar{\odot}}$ whose Trotter product gives back $\xi^\odot$.  As a byproduct existence of projection morphisms onto the subsystems shows that the canonical injections ${E^\ell}^\odot\rightarrow\bfam{{E^1}\circledcirc{E^2}}^\odot$ have an adjoint so that both the injections and the projections are proper spatial morphisms in the sense of Definition \ref{spatdef}.

\bprop
There exist projection morphisms ${p^\ell}^\odot$ onto the subsystems ${E^\ell}^\odot$. These projection morphisms are spatial and continuous.
\eprop

\proof
Let $\breve{p}^\ell_t$ denote the projection in $\sB^{a,bil}(\breve{F}_t)$ onto $E^\ell_t$. One easily checks that $(\breve{p}^\ell_s\odot\breve{p}^\ell_t)(\breve{F}_s\odot\breve{F}_t)=E^\ell_{s+t}=\breve{p}^\ell_{s+t}\breve{F}_{s+t}$. Therefore, the mappings $p^\ell_t$ defined by setting
\beqn{
p^\ell_t(x_{t_n}^n\odot\ldots\odot x_{t_1}^1)
~=~
\breve{p}^\ell_{t_n}x_{t_n}^n\odot\ldots\odot\breve{p}^\ell_{t_1}x_{t_1}^1
}\eeqn
define a projection morphism ${p^\ell}^\odot=\bfam{p^\ell_t}_{t\in\R_+}$ onto ${E^\ell}^\odot\subset \bfam{{E^1}\circledcirc{E^2}}^\odot$. (If $p^\ell_t$ are well-defined then they clearly form a morphism. We refer the reader to [\refcite{BhSk00}, Appendix A] or [\refcite{Ske01}, Appendix A.10] for details about how to well-define mappings on inductive limits.) By [\refcite{BBLS04}, Lemma 5.3.1] (or [\refcite{Ske01}, Lemma 13.2.6]) this morphism is \hl{continuous}, i.e.\ it sends units in $\sU_\om(\bfam{{E^1}\circledcirc{E^2}}^{\sbar{\odot}})$ to units in $\sU_\om({E^\ell}^{\sbar{\odot}})$.~\qed

\blem \label{plem}
The net
\beqn{
(p^1_{t_n}+p^2_{t_n}-\om_{t_n}\om_{t_n}^*)\odot\ldots\odot(p^1_{t_1}+p^2_{t_1}-\om_{t_1}\om_{t_1}^*)
}\eeqn
converges strongly over $\et=(t_n,\ldots,t_1)\in\bJ_t$ to $\id_{(E^1\circledcirc E^2)_t}$.
\elem

\proof
It is not difficult to check that $(p^1_t-\om_t\om_t^*),(p^2_t-\om_t\om_t^*),\om_t\om_t^*$ is a triple of orthogonal projections. Therefore, $p^1_t+p^2_t-\om_t\om_t^*=(p^1_t-\om_t\om_t^*)+(p^2_t-\om_t\om_t^*)+\om_t\om_t^*$ is a projection so that also the net consists entirely of projections. In particular, the net is bounded and it is sufficient to check strong convergence on the total subset of vectors $x_t$ of the form \eqref{genel}. On these vectors convergence is clear, because as soon as the partition $\et$ of the net is finer than the partition of $x_t$ the elements of the net act as identity on $x_t$.~\qed

\bthm \label{uprod}
$\sU_\om(\bfam{{E^1}\circledcirc{E^2}}^{\sbar{\odot}})$ consists of all units $\xi^\odot=\bfam{\xi^1\circledcirc\xi^2\circledcirc\om^{-\beta_\xi}}^\odot$ where ${\xi^\ell}^\odot:=p^\ell\xi^\odot$ $(\ell=1,2)$. Moreover, when $\xi^\odot$ is exponential ${\xi^\ell}^\odot$ are the unique exponential units for ${E^\ell}^{\sbar{\odot}}$ fulfilling $\xi^\odot=\bfam{\xi^1\circledcirc\xi^2}^\odot$.
\ethm

\proof
By Lemmata \ref{plem} and \ref{netlem} we have
\baln{
\xi_t
&~=~
\lim_\et\bfam{(p^1_{t_n}+p^2_{t_n}-\om_{t_n}\om_{t_n}^*)\odot\ldots\odot(p^1_{t_1}+p^2_{t_1}-\om_{t_1}\om_{t_1}^*)}\xi_t
\\
&~=~
\lim_\et(\xi^1_{t_n}+\xi^2_{t_n}-\om_{t_n}e^{t_n\beta_\xi})\odot\ldots\odot(\xi^1_{t_1}+\xi^2_{t_1}-\om_{t_1}e^{t_1\beta_\xi})
\\
&~=~
(\xi^1\bplus\xi^2\bplus-\om^{\beta_\xi})_t.
}\eeqn
As short look at the generator shows that $\bfam{\xi^1\bplus\xi^2\bplus-\om^{\beta_\xi}}^\odot=\bfam{\xi^1\circledcirc\xi^2\circledcirc\om^{-\beta_\xi}}^\odot$.
If $\xi^\odot$ is exponential so are ${\xi^\ell}^\odot$ because ${p^\ell}^\odot$ are spatial and, therefore, $\beta_{\xi^\ell}=\beta_\xi=0$. As ${p^\ell}^\odot$ gives ${\xi^\ell}^\odot$ such exponential units are unique.~\qed

\brem
It is clear that the construction of this section may be generalized to an arbitrary number (finite or infinite) of spatial product systems and that it is associative and commutative (up to canonical isomorphism).
\erem

We close this section by showing that our product is a coproduct rather than a product in the cateogory of spatial product systems with possibly unbounded and possibly not adjointable morphisms that respect the reference unit.

\bthm \label{coprodthm}
Let $({E^\ell}^\odot,{\om^\ell}^\odot)$ $(\ell=1,2)$ and $(E^\odot,\om^\odot)$ denote spatial product systems and let ${w^\ell}^\odot\colon E^\odot\rightarrow{E^\ell}^\odot$ be (possibly unbounded and possibly not adjointable) spatial morphisms. Then with the canonical embeddings ${j^\ell}^\odot\colon{E^\ell}^\odot\rightarrow\bfam{{E^1}\circledcirc{E^2}}^\odot$ there exists a unique (possibly unbounded and possibly not adjointable) spatial morphism
\beqn{
w^\odot
\colon
\bfam{{E^1}\circledcirc{E^2}}^\odot
~\longrightarrow~
E^\odot
}\eeqn
such that $w{j^\ell}^\odot={w^\ell}^\odot$.
\ethm

\proof
Obviously, $w^\odot$ is determined uniquely. Let $\breve{w}_t$ denote the bilinear operator from $\breve{F}_t$ to $E_t$ that sends $x\in E_t^\ell\subset\breve{F}_t$ to $w^\ell_tx$. On elements of the form \eqref{genel} we put $w_t(x^n_{t_n}\odot\ldots\odot x^1_{t_1})=\breve{w}_{t_n}x^n_{t_n}\odot\ldots\odot\breve{w}_{t_1}x^1_{t_1}$. That (well-)defines a bilinear operator $w_t$ from $\bfam{{E^1}\circledcirc{E^2}}_t$ to $E_t$. Clearly, the $w_t$ form a (possibly unbounded and possibly not adjointable) spatial morphism that fulfills the requirments.~\qed

\brem
Even if each ${w^\ell}^\odot$ is bounded, then $w^\odot$ need not be bounded. (For instance, if ${E^\ell}^\odot=E^\odot$ and ${w^\ell}^\odot=\id_{E^\odot}^\odot$ then $w^\odot$ is bounded only in rare occasions. A sufficient criterion, that often is also necessary, is that every $\breve{w_t}$ be a contraction.) Therefore, $\bfam{{E^1}\circledcirc{E^2}}^{\sbar{\odot}}$ is not the coproduct in the category of spatial product systems of Hilbert \nbd{\cB}\nbd{\cB}modules.
\erem

\brem
$w_t$ is the inductive limit of adjointable mappings, but inductive limits of adjointable mappings need not be adjointable. This remains true even if $w_t$ is bounded and extends, thus, to the norm completion of the inductive limit. Only for von Neumann modules we know that every bounded operator has an adjoint. But it remains the fact discussed in the preceding remark that $w_t$ need not be bounded. However, if all $w_t$ have an adjoint then $w^*_t$ is the morphism that would be required by the universal property of a product for the two morphisms ${{w^\ell}^*}^\odot$ within the category of spatial product systems.
\erem

Knowing all this we will appreciate better the product and coproduct properties of completely spatial product systems which we will discuss in Theorem \ref{univp}

\section{Time ordered Fock module}\label{toSEC}

In this section we discuss the time ordered Fock module and show that every completely spatial product system is a time ordered system. Key ingredients are that the generator of the CPD-semigroup associated with a time ordered system has CE-form and the Kolomogorov decomposition of the CPD-part $\eL_0$ of the generator of a completely spatial product systems in \eqref{msgen}. This enables us to define the \it{index} of a spatial product system. Then we use the results from Section \ref{tpuSEC} to show that the index is additive under our product.

Let $F$ be a Hilbert \nbd{\cB}\nbd{\cB}module. By $L^2(\R_+,F)$ we denote the completion of the exterior tensor product $F\otimes L^2(\R_+)$ (and similarly for other measure spaces). $L^2(\R_+,F)$ is a Hilbert \nbd{\cB}\nbd{\cB}module with obvious structure; see [\refcite{Ske01}] for details. As usual, we have $L^2(\R_+,F)^{\sbar{\odot}n}$ $=L^2(\R_+^n,F^{\sbar{\odot}n})$.

By $\Delta_n$ we denote the indicator function of the subset $\bCB{(t_n,\ldots,t_1)\colon t_n>\ldots>t_1>0}$ of $\R_+^n$. Clearly, $\Delta_n$ acts as a projection on $L^2(\R_+,F)^{\sbar{\odot}n}$. The \hl{time ordered Fock module} is
\beq{ \label{todef}
\DG(F)
~=~
\coplus_{n\in\N_0}\Delta_nL^2(\R_+,F)^{\sbar{\odot}n}
}\eeq
where $L^2(\R_+,F)^{\sbar{\odot}0}=\cB$ and $\om=\U\in\cB=L^2(\R_+,F)^{\sbar{\odot}0}$ is the \hl{vacuum}. Denote by $\DG_t(F)$ the restriction of $\DG(F)$ to $\RO{0,t}$. By [\refcite{BhSk00}] $\DG^\odot(F)=\bfam{\DG_t(F)}_{t\in\R_+}$ is a product system of Hilbert modules, the \hl{time ordered system}. The isomorhpism $\DG_s(F)\sbar{\odot}\DG_t(F)\rightarrow\DG_{s+t}(F)$ is obtained by first shifting the interval $\RO{0,s}$ of the left factor to $\RO{t,t+s}$ and then taking the pointwise tensor product of module-valued functions; see [\refcite{Ske01}, Theorem 7.1.3]:
\beqn{
\SB{F_s\odot G_t}(s_m,\ldots,s_1,t_n,\ldots,t_1)
~=~
F_s(s_m-t,\ldots,s_1-t)\odot G_t(t_n,\ldots,t_1).
}\eeqn

\brem
Also the algebraic time ordered Fock module, where in \eqref{todef} direct sum and tensor products are algebraic, gives rise to an (algebraic) product system; see [\refcite{Ske01}]. However, in this product system there are usually (that is, unless $E^{\odot n}=\zero$ for some $n$) no units with components outside the vacuum. Later on we will see another algebraic subsystem of $\DG^\odot(F)$ consisting of the spaces $\DG^{\sU_c}_t(F)$ that are generated algebraically by the continuous units.
\erem

$\DG^\odot(F)$ has a central unital unit $\om^\odot=\bfam{\om_t}_{t\in\R_+}$ with $\om_t=\om$ and we will think always of the time ordered system as the spatial product system $(\DG^\odot(F),\om^\odot)$. But $\DG^\odot(F)$ has lots of other units. Let $\beta\in\cB$ and $\zeta\in F$. Let $\xi^0_t=e^{t\beta}\in\cB=L^2(\R_+,F)^{\sbar{\odot}0}$ be the semigroup in $\cB$ with generator $\beta\in\cB$ and define $\xi^n_t\in\Delta_nL^2(\R_+,F)^{\sbar{\odot}n}$ by setting
\beqn{
\xi^n_t(t_n,\ldots,t_1)
~=~
\xi^0_{t-t_n}\zeta\odot\xi^0_{t_n-t_{n-1}}\zeta\odot\ldots\odot\xi^0_{t_2-t_1}\zeta\xi^0_{t_1}.
}\eeqn
By [\refcite{LiSk01}] $\xi^\odot(\beta,\zeta)=\bfam{\xi_t(\beta,\zeta)}_{t\in\R_+}$ with $\xi^n_t$ being the component of $\xi_t(\beta,\zeta)$ in the \nbd{n}particle sector, defines a unit for $\DG^\odot(F)$. Other results from [\refcite{LiSk01}] can be rephrased in the following way; see [\refcite{Ske01}, Chapter 7].

\bthm\label{LSthm}
$\sU_c(F):=\bCB{\xi^\odot(\beta,\zeta)\colon\beta\in\cB,\zeta\in F}$ is equal to $\sU_\om(\DG^\odot(F))$. Identifying this as the set $\cB\times F$, the generator of $\eU\upharpoonright\sU_c(F)$ is the CE-generator $\eL$ given by
\beqn{
\eL^{(\beta,\zeta),(\beta',\zeta')}(b)
~=~
\AB{\zeta,b\zeta'}+\beta^*b+b\beta'.
}\eeqn
\ethm

We find the form of the maximal completely spatial subsystem of a spatial system and, in particular, of a completely spatial system.

\bthm \label{msssthm}
Let $(E^{\sbar{\odot}},\om^\odot)$ be a spatial product system of Hilbert modules. Then there is a (unique up to two-sided isomorphism) Hilbert \nbd{\cB}\nbd{\cB}module $F$ such that the maximal completely spatial subsystem $({E^{\sU_\om}}^{\sbar{\odot}},\om^\odot)$ of $E^{\sbar{\odot}}$ is isomorphic to $\DG^\odot(F)$. In particular, completely spatial product systems of Hilbert modules are time ordered.
\ethm

\bdefi\label{indexdef}
We refer to the space $F$ as the \hl{index} of a spatial product system.
\edefi

\proof[Proof of Theorem \ref{msssthm}.]
Let $\eL_0$ be the completely positive definite part of the generator $\eL$ of $\eU\upharpoonright\sU_\om(E^{\sbar{\odot}})$ as in \eqref{msgen}. Let $(F,\bfam{\zeta_\xi}_{\xi\in\sU_\om(E^{\sbar{\odot}})})$ be the (completion of the) Kolmogorov decomposition for $\eL_0$, i.e.\ $\eL_0^{\xi,\xi'}(b)=\AB{\zeta_\xi,b\zeta_{\xi'}}$ and $F=\cls{b\zeta_\xi b'}$. Then
\beqn{
\xi^\odot
~\longmapsto~
\xi^\odot(\beta_\xi,\zeta_\xi)
}\eeqn
defines an isometric morphism of the maximal completely spatial subsystem of $E^{\sbar{\odot}}$ onto a subsystem of $\DG^\odot(F)$ (the generators of the associated CPD-semi\-groups coincide and this determins inner products on the total subset of elements of the form \eqref{tensgen}) sending $\om^\odot$ to $\om^\odot$ (obviously, $\zeta_\om=0$).

To see surjectivity, we observe that by Proposition \ref{utoe} the subsystem generated by $\xi^\odot(\beta,\zeta)$ and $\om^\odot$ contains the exponential unit $\xi^\odot(0,\zeta)$ and, more generally, by Lemma \ref{Tlem} any unit $\xi^\odot(0,\lambda\zeta)$ $(\lambda\in\SB{0,1})$. Differentating the continuous function $\lambda\mapsto\xi_t(0,\lambda\zeta)$ with respect to $\lambda$ and putting $\lambda=0$, we obtain the function $\zeta\I_{\RO{0,t}}$ (with $t\in\R_+,\zeta\in F$ arbitrary) in the one-particle sector. By taking tensor products we obtain all functions
\beqn{
\zeta_n\I_{\RO{0,t_n}}\odot\ldots\odot\zeta_1\I_{\RO{0,t_1}}
}\eeqn
$n\in\N,\et\in\bJ_t,\zeta_i\in F$. In the proof of [\refcite{Ske01}, Theorem 7.2.2] we showed that these functions are total in $\DG_t(F)$.

Now let $w^\odot\colon\DG^\odot(F)\rightarrow\DG^\odot(F')$ be a spatial isomorphism between two time ordered product systems. It sends $\sU_c(F)$ to a continuous set of units containing $\om'^\odot$. By Lemma \ref{contlem} this subset is contained in $\sU_c(F')$ so that $w^\odot$ is a (possibly unbounded) continuous morphism. By [\refcite{BBLS04}, Theorem 5.2.1] (or [\refcite{Ske01}, Theorem 13.2.1]) the (possibily unbounded) continuous morphisms (not necessarily spatial) are in one-to-one correspondence with matrices $\G=\smash{\sMatrix{\gamma&\eta^*\\[-.5ex]\eta'&a}}\in\sB^{a,bil}(\cB\oplus F,\cB\oplus F')$. (The morphism $w^\odot$ with matrix $\G$ acts on units as $w\xi^\odot(\beta,\zeta)=\xi^\odot\bfam{\,\gamma+\beta+\AB{\eta,\zeta}\,,\,\eta'+a\zeta\,}$.) By [\refcite{BBLS04}, Corollary 5.2.4] (or [\refcite{Ske01}, Corollary 13.2.4]) in order that $w^\odot$ be a (not necessarily spatial) isomorphism, $a\in\sB^{a,bil}(F,F')$ must be a two-sided unitary, i.e.\ an isomorphism $F\rightarrow F'$.~\qed

\brem
We check immediately that in order that $w^\odot$ be spatial, we must have $\gamma=0,\eta=0,\eta'=0$.
\erem

We use the concrete form of the morphisms in order to show that the spatial structure of $\DG^\odot(F)$ does not depend on the choice of the reference unit.

\bprop\label{transcsrem}
If  $\xi^\odot=\xi^\odot(\beta,\zeta)$ is another central unital unit then there is an automorphism of $\DG^\odot(F)$ sending $\om^\odot$ to $\xi^\odot$.
\eprop

\proof
In order that $\xi^\odot$ be central, it is necessary and sufficient that $\beta\in C_\cB(\cB),\zeta\in C_\cB(F)$. (This follows by investigating carefully the generator in Theorem \ref{LSthm} taking into account that $\eL^{(\beta',\zeta'),(\beta,\zeta)}(b)=\eL^{(\beta',\zeta'),(\beta,\zeta)}(\U)b$ for all $b,\beta'\in\cB;\zeta'\in F$.) In order that $\xi^\odot$ be unital it is necessary and sufficient that
\beqn{
0
~=~
\eL^{\xi,\xi}(\U)
~=~
\AB{\zeta,\zeta}+\beta^*+\beta.
}\eeqn
This means the real part of $\beta$ is $-\frac{\AB{\zeta,\zeta}}{2}$ and the imaginary part is arbitrary. The conditions in [\refcite{BBLS04}, Corollary 5.2.4] that the endomorphism $w^\odot$ determined by the matrix $\G=\smash{\sMatrix{\gamma&\eta^*\\[-.5ex]\eta'&a}}\in\sB^{a,bil}(\cB\oplus F)$ be an automorphism are that $a$ be an arbitrary automorphism of $F$ and $\eta'$ be an arbitrary element of $C_\cB(F)$ while $\eta=-a^*\eta'$ and $\gamma=ih-\frac{\AB{\eta',\eta'}}{2}$ for some arbitrary self-adjoint element $h\in C_\cB(\cB)$. Clearly, $w^\odot$ sends $\om^\odot=\xi^\odot(0,0)$ to $\xi^\odot(\gamma,\eta')$. Therefore, if we choose $\eta'=\zeta$ and $\gamma=\beta$ ($a$ arbitrary and $\eta$ correspondingly), then we obtain an automorphism sending $\om^\odot$ to $\xi^\odot$.~\qed

\lf
Now we show that $F$ merits to be called an index.

\bthm\label{indsumthm}
Let $({E^\ell}^{\sbar{\odot}},{\om^\ell}^{\sbar{\odot}})$ $(\ell=1,2)$ be two spatial product systems with indices $F^\ell$. Then the index of $(\bfam{{E^1}\circledcirc{E^2}}^{\sbar{\odot}},\om^\odot)$ is $F^1\oplus F^2$. In particular, $(\bfam{\DG(F^1)\circledcirc\DG(F^2)}^{\sbar{\odot}},\om^\odot)$ is isomorphic to $(\DG^\odot(F^1\oplus F^2),\om^\odot)$.
\ethm

\proof
By Theorem \ref{uprod} any continuous unit for $\bfam{{E^1}\circledcirc{E^2}}^{\sbar{\odot}}$ may be obtained as a Trotter product of units in ${E^\ell}^{\sbar{\odot}}\subset\bfam{{E^1}\circledcirc{E^2}}^{\sbar{\odot}}$, and by Theorem \ref{msssthm} the units are even from the maximal completely spatial subsystems $\DG^\odot(F^\ell)\subset{E^\ell}^{\sbar{\odot}}$.By Corollary \ref{etou} we may restrict to exponential units. By looking at the generators of the relevant CPD-semigroups we see that sending the Trotter product of exponential units $\xi^\odot(0,\zeta^\ell)$ for $\DG^\odot(F^\ell)$ to the exponential unit $\xi^\odot(0,\zeta^1\oplus\zeta^2)$ for $\DG^\odot(F^1\oplus F^2)$ defines a surjective isometric morphism (i.e.\ an isomorphism) from the maximal completely spatial subsystem $\bfam{\DG(F^1)\circledcirc\DG(F^2)}^{\sbar{\odot}}$ of $\bfam{{E^1}\circledcirc{E^2}}^{\sbar{\odot}}$ onto $\DG^\odot(F^1\oplus F^2)$. Clearly, the reference units are preserved.~\qed

\lf
In the case of type I systems of Hilbert spaces our product is nothing but the tensor product. In the case of type II systems we obtain at least a subsystem of the tensor product. (It may coincide with the tensor product, but it need not.) To see this we use for $x_t\in E^1_t\cup E^2_t\subset E^1_t\oplus E^2_t$ the notation
\beqn{
(x_t)^\ell
~=~
\begin{cases}
x_t&\text{~for $x_t\in E^\ell_t$}\\
\om^\ell_t&\text{~otherwise}.
\end{cases}
}\eeqn

\bprop\label{ntprem}
Let $({\eH^\ell}^{\sbar{\otimes}},{\om^\ell}^\otimes)$ $(\ell=1,2)$ be two spatial product systems of Hilbert spaces and denote by $(F^{\sbar{\otimes}},\om^\otimes)$ their product. Then the mapping
\beqn{
u_t
\colon
x_{t_n}^n\otimes\ldots\otimes x_{t_1}^1
~\longmapsto~
((x_{t_n}^n)^1\otimes\ldots\otimes(x_{t_1}^1)^1)\otimes((x_{t_n}^n)^2\otimes\ldots\otimes(x_{t_1}^1)^2)
}\eeqn
extends as an isometry $F_t\rightarrow\eH^1_t{\sbar{\otimes}}\eH^2_t$ and the family $u^\otimes=\bfam{u_t}_{t\in\R_+}$ is an isometric morphism of product systems. In the case of type I systems it is an isomorphism, i.e.\ the $u_t$ are unitary.\footnote{When the spatial product systems of Hilbert spaces are continuous in the sense of [\refcite{Ske03b}], then one may show that already one completely spatial factor is sufficient to have equality of our product with the tensor product. In the separable case (i.e.\ Arveson systems) this can be shown using results from Liebscher [\refcite{Lie00p1}]. (Treating the general case would require to repeat a good deal of [\refcite{Ske03b}, Section 7], so we decided not to include a proof here.) Liebschers methods also show that there exist examples when both spatial factors are not completely spatial where our product is a proper subsystem of the tensor product. (The two subsystems $\bfam{\eH^1_t\otimes\om^1_t}_{t\in\R_+}$ and $\bfam{\om^2_t\otimes\eH^2_t}_{t\in\R_+}$ do not generate the whole tensor product $\bfam{\eH^1\sbar{\otimes}\eH^2}^{\sbar{\otimes}}$.) A proof of this statement will appear in Bhat, Liebscher and Skeide [\refcite{BLS02p}].}
\eprop

\proof
The mapping is isometric, because the mapping $x_t\mapsto(x_t)^1\otimes(x_t)^2$ on $\breve{F}_t$ is isometric. In the case of type I systems $\DG^\otimes(K^\ell)$ the range of $u_t$ contains all exponential vectors in $\DG_t(K^1)\sbar{\otimes}\DG_t(K^2)=\DG_t(K^1\oplus K^2)$ to step functions with values in the total set $K^1\cup K^2\ni 0$ which are total by [\refcite{Ske01}, Theorem 7.4.3]. (See also Parthasarathy and Sunder [\refcite{PaSu98}] and Skeide [\refcite{Ske00a}].)~\qed

\lf
We close this section by showing that the product merits to be called a product (actually also a coproduct) at least in the category of completely spatial product systems. In order to avoid problems with unbounded morphisms we consider spatial morphisms of the algebraic subsystems ${\DG^{\sU_c}}^\odot(F)$. From [\refcite{BBLS04}] (see [\refcite{Ske01}, Theorem 13.2.1]) it follows that such morphisms $w^\odot\colon{\DG^{\sU_c}}^\odot(F^1)\rightarrow{\DG^{\sU_c}}^\odot(F^2)$ are in one-to-one correspondence with operators $a\in\sB^{a,bil}(F^1,F^2)$, and that $w^\odot$ acts on a unit as $w\xi^\odot(\beta,\zeta)=\xi^\odot(\beta,a\zeta)$. If $a$ is a contraction, then $w^\odot$ is bounded (even contractive). The converse need not be true.

\bthm \label{univp}
Let $\DG^\odot(F^\ell)$ $(\ell=1,2)$ and $\DG^\odot(F)$ be completely spatial product systems and let ${w^\ell}^\odot\colon{\DG^{\sU_c}}^\odot(F)\rightarrow{\DG^{\sU_c}}^\odot(F^\ell)$ be (possibly unbounded) morphisms. Then
\begin{enumerate}
\item
$\DG^\odot(F^1\oplus F^2)$ with the projection morphisms ${p^\ell}^\odot$ is the product in the category of completely spatial product systems, i.e.\ there exists a unique (possibly unbounded) morphism
\beqn{
w^\odot
\colon
{\DG^{\sU_c}}^\odot(F)
~\longrightarrow~
{\DG^{\sU_c}}^\odot(F^1\oplus F^2)
}\eeqn
such that $p^\ell w^\odot={w^\ell}^\odot$.

\item
$\DG^\odot(F^1\oplus F^2)$ with the canonical embeddings ${j^\ell}^\odot$ is the coproduct in the category of completely spatial product systems, i.e.\ there exists a unique (possibly unbounded) morphism
\beqn{
{w^*}^\odot
\colon
{\DG^{\sU_c}}^\odot(F^1\oplus F^2)
~\longrightarrow~
{\DG^{\sU_c}}^\odot(F)
}\eeqn
such that $w^*{j^\ell}^\odot={w^{\ell^*}}^\odot$.
\end{enumerate}
${w^*}^\odot$ is, indeed, the adjoint of $w^\odot$.
\ethm

\proof
Let $a^\ell\in\sB^{a,bil}(F,F^\ell)$ be the operators generating the morphisms ${w^\ell}^\odot$. It follows that the operator $a=a^1+a^2\in\sB^{a,bil}(F,F^1\oplus F^2)$ generates a morphism $w^\odot$ which has all the properties desired for the product, and that its adjoint (generated by $a^*$) has all the properties desired for the coproduct.~\qed

\brem
Even if $a^\ell$ are contractions, then $a=a^1+a^2$ need not be. Nevertheless, $\DG^\odot(F^1\oplus F^2)$ is determined by each of the preceding universal properties up to isomorphism. This is so, because in the isomorphism (constructed in the usual way from the univeral property) between two candidates the operator $a$ decomposes as $a^1\oplus a^2$ where $a^\ell$ are unitarily equivalent to unitaries in $\sB^a(F^\ell)$. (This is so, because the restriction to $\DG^\odot(F^\ell)$ must define an isomorphism of $\DG^\odot(F^\ell)$.) Therefore, also $a$ is a unitary.
\erem

\section{Noises and spatial product systems}\label{wnspsSEC}

{\parskip0.5ex plus 0.5ex minus 0.5ex
\nbd{E_0}semigroups acting as strict unital endomorphisms on $\sB^a(E)$ for some Hilbert \nbd{\cB}module $E$ give rise to a product system $E^{\sbar{\odot}}$ of Hilbert \nbd{\cB}modules. In the case when $E=H$ is a separable infinite-dimensional Hilbert space the correspondence between strongly continuous \nbd{E_0}semigroups (up to cocycle conjugacy) and Arveson systems (up to isomorphism) is one-to-one. For Hilbert modules one direction of this result is wrong. (There are \nbd{E_0}semigroups on $E^1$ and $E^2$ having the same product system but where $\sB^a(E^1)$ and $\sB^a(E^2)$ are nonisomorphic so that cocycle conjugacy has no meaning. Only if we fix the isomorphism class of $E$, then two \nbd{E_0}semigroup have the same product system, if and only if they are concycle conjugate; see [\refcite{Ske02}].\footnote{If we fix, instead, the strict isomorphism class of $\sB^a(E)$, then the product system of two conjugate and therefore also of cocycle conjugate \nbd{E_0}semigroups are Morita equivalent; see Skeide [\refcite{Ske04p}].}) The other direction, constructing an \nbd{E_0}semigroup from a product system, promisses to remain true (if we do not ask too ingenuously), but presently we do not yet have a proof. This difficulty dissappears, however, as soon as we have a (unital or continuous) unit.

In this section we repeat results concerning \nbd{E_0}semigroups on Hilbert modules with particular emphasis on the spatial case that corresponds to noises.\footnote{The term \it{noise} is justified by the fact that such \nbd{E_0}semigroups come along with filtrations of subalgebras that are \hl{monotone independent} in a certain invariant conditional expectation; see Skeide [\refcite{Ske04}].} In the following section we will apply the results to extend the product of spatial product systems to a product of noises.\footnote{We discuss here the construction from Skeide [\refcite{Ske02}] of a product system  from an \nbd{E_0}semigroup on $\sB^a(E)$ based on the assumption that $E$ has a unit vector $\xi$. (This construction generalizes directly the construction of Arveson systems in Bhat [\refcite{Bha96}]. It has nothing to do with Arveson's original construction from [\refcite{Arv89}]. In fact, we explain in Skeide [\refcite{Ske03b,Ske03p1}] that Arveson's construction leads to product systems that are anti-isomorphic to that from Bhat's construction, and we discuss a generalization of Arveson's construction that works only for von Neumann modules and leads to product systems of bimodules over the commutant $\cB'$ of $\cB$.) Only recently, we have freed the construction from assuming existence of a unit vector (and that $\cB$ be unital) in Muhly, Solel and Skeide [\refcite{MSS03p}]. The old approach from [\refcite{Ske02}] has the advantage of very concrete identifications of the members of the product system as submodules of $E$. (Additionally, requiring for spatial product systems existence of intertwining semigroups of isometries in the spirit of Powers, as a consequence, will lead to unital algebras and loads of unit vectors.)}

In [\refcite{Ske02}] (see [\refcite{Ske01}, Section 14.1]) we associated with a triple $(E,\vt,\xi)$ consisting of a Hilbert \nbd{\cB}module, a \hl{strict} (i.e.\ $\vt_t$ is \nbd{*}strongly continuous on bounded subsets of $\sB^a(E)$ for all $t\in\R_+$) \nbd{E_0}semigroup $\vt$ on $\sB^a(E)$ and a unit vector $\xi\in E$ a product system $E^{\sbar{\odot}}$ in the following way. We define a family $j=\bfam{j_t}_{t\in\R_+}$ of representations $j_t$ of $\cB$ on $E$ by setting $j_0(b)=\xi b\xi^*$ and $j_t=\vt_t\circ j_0$. Then with $p_t=j_t(\U)$ the submodule $E_t=p_tE$ of $E$ becomes a two-sided Hilbert module with left action $b.x_t=j_t(b)x_t$. Then one checks that
\beq{\label{srf}
x\odot y_t
~=~
\vt_t(x\xi^*)y_t
}\eeq
defines a unitary identification $E\sbar{\odot}E_t=E$ such that $a\odot\id_{E_t}=\vt(a)$. (One easily verifies that the mapping $x\odot y_t\mapsto\vt_t(x\xi^*)y_t$ is isometric. Surjectivity is slightly harder to verify and uses that $\vt_t$ is strict.) Restricting \eqref{srf} to the subspace $E_s\sbar{\odot}E_t$ of $E\sbar{\odot}E_t$, we obtain two-sided isomorphisms $E_s\sbar{\odot}E_t=E_{s+t}$. By the semigroup property of $\vt$ we see that
\beqn{
(E\sbar{\odot}E_s)\sbar{\odot}E_t
~=~
E\sbar{\odot}(E_s\sbar{\odot}E_t).
}\eeqn
Restricting to $E_r\subset E$ we see that $(E_r\sbar{\odot}E_s)\sbar{\odot}E_t=E_r\sbar{\odot}(E_s\sbar{\odot}E_t)$ so that $E^{\sbar{\odot}}=\bfam{E_t}_{t\in\R_+}$ is a product system; see [\refcite{Ske01}, Theorem 14.1.1].

$(E,\vt,\xi)$ is a \hl{weak dilation} (i.e.\ the mappings $T_t\colon b\mapsto\AB{\xi,j_t(b)\xi}$ form a semigroup $T$ on $\cB$, necessarily completely positive and unital), if and only if the $p_t$ are increasing (i.e.\ $p_t\ge p_0$ for all $t\in\R_+$). In this case the vectors $\xi_t=\xi$ are in $E_t$ for all $t\in\R_+$ and $\xi^\odot=\bfam{\xi_t}_{t\in\R_+}$ is a unit for $E^{\sbar{\odot}}$ that satisfies also $\xi=\xi\odot\xi_t$; see [\refcite{Ske01}, Proposition 14.1.6]. We deonote $E_\infty=\ol{\bigcup_{t\in\R_+}E_t}$. The weak dilation $(E,\vt,\xi)$ is called \hl{primary}, if $E=E_\infty$. In this case, $\lim_{t\to\infty}p_t=\id_E$ strongly. The weak dilation $(E,\vt,\xi)$ is a \hl{(weak) noise}, i.e.\ $T$ is the trivial semigroup on $\cB$, if and only if $\xi^\odot$ is also central. In this case, $\vp_0(a)=p_0bp_0$ defines a conditional expectation $\vp_0\colon\sB^a(E)\rightarrow j_0(\cB)$ that is \hl{invariant} for $\vt$, i.e.\ $\vp_0\circ\vt_t=\vp_0$ for all $t\in\R_+$. Observe that if $(E,\vt,\xi)$ is primary then the strong limit $\lim_{t\to\infty}j_t(b)$ defines a unital representation on $E=E_\infty$ so that $E$ is turned into a Hilbert \nbd{\cB}\nbd{\cB}module. In other words, $\cB$ is identified as a unital subalgebra of $\sB^a(E)$ which, clearly, is invariant for $\vt$. Since $j_t(b)\xi=\xi b$ for all $t$, we find that $b\xi=\xi b$. In other words, $\vp(a)=\AB{\xi,a\xi}$ defines a \hl{vector expectation} (in analogy with \hl{vector state}) onto $\cB\subset\sB^a(E)$. (In fact, $(E,\xi)$ may be considered as the \hl{GNS-construction} of the conditional expectation $\vp$ in the sense of [\refcite{Pas73}].) Since $j_0\circ\vp=\vp_0$ and $\vp\circ\vp_0=\vp$, the expectation is invariant for $\vt$. In general, we call $(E,\vt,\xi)$ a \hl{unital noise}, if $E$ is a Hilbert \nbd{\cB}\nbd{\cB}module and $\vp=\AB{\xi,\bullet\xi}$ is an invariant vector expectation onto the subalgebra $\cB$ of $\sB^a(E)$.

\lf
So far, we have investigated how to construct product systems from \nbd{E_0}semigroups, and how additional structure like being a dilation or even a noise gives rise to additional structure of the product system like possessing a unit or even being spatial. Now we discuss in how far we may reverse the construction. As we do not yet know, whether every (reasonable) unitless product system comes from an \nbd{E_0}semigroup, we discuss only the (rather simple) case in presence of a (unital) unit.

From [\refcite{BhSk00}] (see [\refcite{Ske01}, Section 11.4]) we know how to construct from a unital unit $\xi^\odot$ in a product system $E^{\sbar{\odot}}$ a Hilbert module $E_\infty=\ol{\limind E_t}$, where $E_t$ is embedded isometrically into $E_{s+t}$ via $x_t\mapsto\xi_s\odot x_t$. By obvious extension of \eqref{tpsid} and its associativity, we find
\beq{ \label{Efacs}
E_\infty
~=~
E_\infty\sbar{\odot}E_t
\text{~~~~~~and~~~~~~}
(E_\infty\sbar{\odot}E_s)\sbar{\odot}E_t
~=~
E_\infty\sbar{\odot}(E_s\sbar{\odot}E_t).
}\eeq
It follows that $\vt_t\colon a\mapsto a\odot\id_{E_t}\in\sB^a(E_\infty\sbar{\odot}E_t)=\sB^a(E_\infty)$ defines a strict \nbd{E_0}semigroup $\vt=\bfam{\vt_t}_{t\in\R_+}$ on $\sB^a(E_\infty)$. Under the inductive limit the vectors $\xi_t\in E_t$ are all mapped to the same unit vector $\xi\in E_\infty$. This shows clearly, that $E_\infty$ is not a two-sided module (giving back the correct left action of $\cB$ on $E_t$), unless $\xi^\odot$ is central. Embedding $\cB$ via $j_0(b)=\xi b\xi^*$ into $\sB^a(E_\infty)$, we find that $\AB{\xi,\vt_t\circ j_0(b)\xi}=\AB{\xi_t,b\xi_t}$, i.e.\ $(E_\infty,\vt,\xi)$ is a weak dilation of the CP-semigroup $T_t=\AB{\xi_t,b\xi_t}$. The dilation is primary and constructing its product system $\bfam{p_tE_\infty}_{t\in\R_+}$ gives back the product system we started with. The dilation is a (weak) noise, if and only if $\xi^\odot$ is also central. In this case, as before $E_\infty$ may be turned into a two-sided Hilbert module such that the $E_t$ are contained as two-sided submodules. In particular, with any spatial product system $(E^{\sbar{\odot}},\om^\odot)$ we may associate a strict primary weak noise $(E_\infty,\vt,\om)$.
}

\section{The product of noises} \label{wnSEC}

So far we have introduced index and product of spatial product systems. But often product systems stem from \nbd{E_0}semigroups, and so far we do not have a product of \nbd{E_0}semigroups in the module case that could play the role of the tensor product of \nbd{E_0}semigroups in the Hilbert space case. In this section we present a product at least of those \nbd{E_0}semigroups that are \it{(weak) noises} as discussed in the preceding section. We will see that if the noises are primary, then also the product will be primary and, therefore, we find also a product of unital noises coming from primary noises. We close the section with the technical result that our product preserves continuity in time.

\bthm\label{WNprodthm}
Let $(E^\ell,\vt^\ell,\om^\ell)$ $(\ell=1,2)$ be two strict weak noises with associated spatial product systems $({E^\ell}^{\sbar{\odot}},{\om^\ell}^\odot)$ and product $(\bfam{E^1\circledcirc E^2}^{\sbar{\odot}},\om^\odot)$. Then there exists a strict weak noise $(F,\vt,\om)$ fulfilling the following properties:
\begin{enumerate}
\item\label{1}
The spatial product system associated with $(F,\vt,\om)$ is $(\bfam{E^1\circledcirc E^2}^{\sbar{\odot}},$ $\om^\odot)$.

\item
$F$ contains $E^\ell$ as submodules in such a way that $\AB{x^1,x^2}=\AB{x^1,\om}\AB{\om,x^2}$ for all $x^1\in E^1,x^2\in E^2$.

\item
$F$ is generated by its submodules $E^\ell$ and the product system $\bfam{E^1\circledcirc E^2}^{\sbar{\odot}}$ in the sense that elements of the form $x^\ell\odot y_t$ $(\ell=1,2;x^\ell\in E^\ell;t\in\R_+;y_t\in(E^1\circledcirc E^2)_t)$ are total in $F$.
\end{enumerate}
Moreover, if $(F',\vt',\om')$ is another strict weak noise fulfilling these properties, then $(F',\vt',\om')$ is unitarily equivalent to $(F,\vt,\om)$, that is there is unitary $u\in\sB^a(F,F')$ intertwining $\vt'$ and $\vt$ and $u\om=\om'$.
\ethm

\bdefi
We call $(F,\vt,\om)$ the \hl{product} of $(E^1,\vt^1,\om^1)$ and $(E^2,\vt^2,\om^2)$ and we denote the product by $(E^1\circledcirc E^2,\vt^1\circledcirc\vt^2,\om)$.
\edefi

\proof[Proof of Theorem \ref{WNprodthm}.~]
Similar to Section \ref{prodSEC} we set $\wh{E}^\ell=E^\ell\ominus\om^\ell\cB$ and $\breve{G}=\om\cB\oplus\wh{E}^1\oplus\wh{E}^2$ with obvious identifications of $E^\ell$ as submodules. We define $G_t=\breve{G}\sbar{\odot}\bfam{E^1\circledcirc E^2}_t$. We observe that $\breve{F}_t=\om_t\cB\oplus\wh{E}^1_t\oplus\wh{E}^2_t\subset\bfam{E^1\circledcirc E^2}_t$ so that $\breve{G}\odot\breve{F}_t\subset\breve{G}\odot\bfam{E^1\circledcirc E^2}_t=G_t$ and as in Equation \eqref{Fst} and its successor we show that $\breve{G}\subset\breve{G}\odot\breve{F}_t$ so that, in the end, $\breve{G}\subset G_t$. Therefore,
\beqn{
G_t
~=~
\breve{G}\sbar{\odot}\bfam{E^1\circledcirc E^2}_t
~\subset~
G_s\sbar{\odot}\bfam{E^1\circledcirc E^2}_t
~=~
\breve{G}\sbar{\odot}\bfam{E^1\circledcirc E^2}_s\sbar{\odot}\bfam{E^1\circledcirc E^2}_t
~=~
G_{s+t}
}\eeqn
and, obviously, the $G_t$ with the canonical embeddings $G_t\rightarrow G_{s+t}$ form an inductive system of Hilbert \nbd{\cB}modules with (completed) inductive limit $G$. Clearly, $G$ with the product system $\bfam{E^1\circledcirc E^2}^{\sbar{\odot}}$ fulfills \eqref{Efacs} and, therefore, by setting $\vt_t(a)=a\odot\id_{(E^1\circledcirc E^2)_t}$ we define a strict \nbd{E_0}semigroup on $\sB^a(G)$. The vector $\om=\om\odot\om_t$ is in all $G_t$ and, thus, gives rise to a unit vector $\om$ in $G$. All three together give rise to a weak noise that fulfills the stated properties.

To see uniqueness we define the unitary $u$ simply by sending element $x^\ell\odot y_t$ from the total subset of $F$ to the corresponding element from the total subset of $F'$. Choosing a rank-one operator $a\in\sB^a(F)$ we see $\vt_t(a)$ and $\vt'_t(uau*)$ act in the same way on elements of the form
\beqn{
x^\ell\odot y_{s_m}\odot\ldots\odot y_{s_1}\odot z_{t_n}\odot\ldots\odot z_{t_1}
}\eeqn
with $t_n+\ldots+t_1=t$, $x^\ell\in E^\ell$ and $y_{s_k}\in E^1_{s_k}$ or $E^2_{s_k}$, $z_{t_j}\in E^1_{t_j}$ or $E^2_{t_j}$. A moments thought shows that these elements still form a total subset. Of course, $u\om=\om'$.~\qed

\bcor
The product is primary, if and only if both factors are. In particular, our product gives rise to a product of (unital) primary weak noises.
\ecor

\proof
$x^\ell\odot y_t$ generate $E^1\circledcirc E^2$. But, if $(E^\ell,\vt^\ell,\om^\ell)$ is primary, then it is sufficient to consider only $x^\ell\in E^\ell_s$. Therefore, $x^\ell\odot y_t\in E^\ell_s\odot(E^1\circledcirc E^2)_t\subset(E^1\circledcirc E^2)_{s+t}$. This shows that $(\vt^1\circledcirc\vt^2)(\om\om^*)$ increases to $\id_{E^1\circledcirc E^2}$.~\qed

\bdefi
The \hl{index} of a strict weak noise is the index of its associated spatial product system.
\edefi

From Property \ref{1} of Theorem \ref{WNprodthm} we obtain immediately:

\bcor
The index of strict weak noises is additive under product.
\ecor

Recall that $\vt$ is \hl{strongly continuous} if $t\mapsto\vt_t(a)x$ is continuous for all $a\in\sB^a(E),x\in E$. (In the case of a von Neumann module $E$ we might think of the \nbd{\sigma}weak topology of the von Neumann algebra $\sB^a(E)$. Here we stay at the level of Hilbert modules and consider only the strongly continuous case.)

We know from [\refcite{Ske01}, Theorem 11.4.12] that the \nbd{E_0}semigroups constructed from a unital unit in a type I system, in particular, those constructed from a completely spatial system, are strongly continuous and, because the product is again completely spatial, also the product of such noises must be strongly continuous. We generalize this to arbitrary strongly continuous noises.

\bthm\label{scprodthm}
Let $(E^\ell,\vt^\ell,\om^\ell)$ $(\ell=1,2)$ be two strict strongly continuous weak noises. Then their product is strongly continuous, too.
\ethm

\proof
The crucial step in the proof [\refcite{Ske01}, Theorem 11.4.12] was the right shift $\s_t\colon x\mapsto x\odot\om_t$. (For an arbitrary unit this mapping is `very bad', in the sense that it is not right linear. Here where we are dealing with central units, $\s_t$ is an isometry which is even adjointable, because there exists the projection $\id\odot\om_t\om_t^*$ onto the range; cf.\ [\refcite{Ske01}, Proposition 1.5.13].)

First, we observe that for each $x\in E^\ell$ separately the mapping $t\mapsto x\odot\om^\ell_t\in E^\ell$ and, therefore, also the mapping $t\mapsto x\odot\om_t\in\bfam{E^1\circledcirc E^2}$, is (norm) continuous. (This follows from $x\odot\om^\ell_t=\vt^\ell_t(x{\om^\ell}^*)\om^\ell$ and strong continuity of $\vt^\ell$.) Applying the projection $p_t$ we find that $x_t\in E^\ell_t$ is close to $p_t(x_t\odot\om^\ell_\ve)=(p_{t-\ve}x_t)\odot\om^\ell_\ve$ which in turn is close to $\om^\ell\odot(p_{t-\ve}x_t)$ for all sufficiently small $\ve\ge0$. Now let
\beq{ \label{total}
X
~=~
x\odot x^n_{t_n}\odot\ldots\odot x^1_{t_1}
}\eeq
be an element in $\bfam{E^1\circledcirc E^2}$ where $x$ and $x^i_{t_i}$ are from $E^1$ or $E^2$ and $E^1_{t_i}$ or $E^2_{t_i}$, respectively, chosen independently. From the preceding considerations it follows that
\baln{
X
~~~\approx~~~~~
x~~~~~\odot(\om_\ve\odot p_{t_n-\ve}x^n_{t_n})\odot\ldots\odot(\om_\ve\odot p_{t_1-\ve}x^1_{t_1})
&
\\
~\approx~
(x\odot\om_\ve)\odot(p_{t_n-\ve}x^n_{t_n}\odot\om_\ve)\odot\ldots\odot(p_{t_1-\ve}x^1_{t_1}\odot\om_\ve)
&
~~~~\approx~
X\odot\om_\ve
}\ealn
for sufficiently small $\ve$. In other words, $\s_t$ is strongly continuous on the total subset of vectors of the form \eqref{total} and, therefore, everywhere. From
\baln{
\norm{(\vt_\ve(a)-a)x}
&~\le~
\norm{\vt_\ve(a)(x-\s_tx)}+\norm{\vt_\ve(a)\s_tx-ax}
\\
&~\le~
\norm{a}\norm{x-\s_tx}+\norm{\s_tax-ax}
~\to~
0
}\ealn
it follows that $\vt$ is strongly continuous.~\qed

\lf
Of course, also the results of this section extend to an arbitrary number of factors, and the constructions are associative and commutative.

\section{The product system of a free flow} \label{freeflow}

\noindent
The time shift endomorphism on $\sB\bfam{\G(L^2(\R_+,K))}$ is also refered to as the \hl{CCR-flow} of index $\dim K$. In analogy we refer to the time shift on $\sB^a(\DG(F))$ as the \hl{(generalized) CCR-flow} of index $F$. Also on the full Fock module $\cF(L^2(\R_+,F))=\coplus_{n\in\N_0}L^2(\R_+,F)^{\sbar{\odot} n}$ we have a time shift $\s_t$ which induces an \nbd{E_0}semigroup $\sS_t$ on $\sB^a\bfam{\cF(L^2(\R_+,F))}$ via the facotrization
\baln{
\cF(L^2(\R_+,F))
&
~=~
\cF(L^2(\RO{t,\infty},F))~\sbar{\odot}~\bfam{\cB\om\,\oplus\,L^2(\RO{0,t},F)\sbar{\odot}\cF(L^2(\R_+,F))}
\\
&
~=~
\cF(L^2(\R_+,F))~\sbar{\odot}~\bfam{\cB\om\,\oplus\,L^2(\RO{0,t},F)\sbar{\odot}\cF(L^2(\R_+,F))}.
}\ealn
Together with the vacuum vector $\om$, the triple $(\cF(L^2(\R_+,F)),\sS,\om)$ is a primary noise which we call the \hl{free flow} of  \hl{free index} $F$.

In this section we show that the associated product system $\cF^\odot(F)=\bfam{\cF_t(F)}_{t\in\R_+}$ is completely spatial and we determine its index. This is a generalization of Fowler [\refcite{Fow95}] who did that program for Hilbert spaces. Like Fowler we construct explicitly all units, observe that they are generating and, after guessing from their form the correct index, we define an explicit isomorphism. It seems, however, that our computations are considerably simpler. On the one hand, our construction of product systems from \nbd{E_0}semigroups works very direct. We just can read off the correct product system from the above factorization by applying the simple identifications as described in Section \ref{wnspsSEC}. On the other hand, we (are forced to) work in the time ordered version of Fock module instead of that of symmetric Fock spaces and it turns out that this simplifies combinatorics considerably. In our opinion this is a strong indication that also in the analysis of Arveson systems it could be more convenient to consider type I Arveson systems as time ordered Fock spaces rather than symmetric Fock spaces.

For some measurable subset $S\subset\R_+$ let us denote $E_S=L^2(S,F)$. We remarked already in [\refcite{Ske01}, Example 14.1.4] that by our construction of product systems from \nbd{E_0}semigroups, $\cF_t(F)=p_t\cF(E_{\R_+})=\cB\om\oplus E_{\RO{0,t}}\sbar{\odot}\cF(E_{\R_+})$. One may check that the identification
\baln{
\bfam{\cB\om
&
\oplus E_{\RO{0,s}}\sbar{\odot}\cF(E_{\R_+})}\sbar{\odot}\bfam{\cB\om\oplus E_{\RO{0,t}}\sbar{\odot}\cF(E_{\R_+})}
\\
&
~\cong~
\s_t\bfam{\cB\om\oplus E_{\RO{0,s}}\sbar{\odot}\cF(E_{\R_+})}\sbar{\odot}\bfam{\cB\om\oplus E_{\RO{0,t}}\sbar{\odot}\cF(E_{\R_+})}
\\
&
~=~
\bfam{\cB\om\oplus E_{\RO{t,t+s}}\sbar{\odot}\cF(E_{[t,\infty)})}\sbar{\odot}\bfam{\cB\om\oplus E_{\RO{0,t}}\sbar{\odot}\cF(E_{\R_+})}
\\
&
~=~
\cB\om\,\oplus\,E_{\RO{0,t}}\sbar{\odot}\cF(E_{\R_+})\,\oplus\,E_{\RO{t,t+s}}\sbar{\odot}\cF(E_{[t,\infty)})\sbar{\odot}\bfam{\cB\om\oplus E_{\RO{0,t}}\sbar{\odot}\cF(E_{\R_+})}
\\
&
~=~
\cB\om\,\oplus\,E_{\RO{0,t}}\sbar{\odot}\cF(E_{\R_+})\,\oplus\,E_{\RO{t,t+s}}\sbar{\odot}\cF(E_{\R_+})
~=~
\cB\om\,\oplus\,E_{\RO{0,t+s}}\sbar{\odot}\cF(E_{\R_+})
}\ealn
gives the correct product system structure. (``$=$'' means canonical identification as subspace of $\cF(E_{\R_+})$.) Of course, the restriction to (separable) Hilbert spaces of this result can be found in [\refcite{Fow95}] (with a different proof).

Our next goal is to find the form of all units. Let $\xi^\odot$ be a unit for $\cF^\odot(F)$ and expand it into $\xi_t=\bigoplus_{n\in\N_0}\xi_t^n$ where $\xi_t^n\in E_{\RO{0,t}}\sbar{\odot}E_{\R_+}^{\sbar{\odot}(n-1)}$. By the unit property we have
\beq{ \label{nunit}
\xi_{s+t}^n
~=~
\sum_{k=0}^n\s_t\xi_s^k\odot\xi_t^{n-k}.
}\eeq
In order to ``derive'' the form of the units, let us start by assuming that the $\xi_t^n$, indeed, are functions of $(t_n,\ldots,t_1)\in\RO{0,t}\times\R_+^{n-1}$. In other words, \eqref{nunit} means equality everywhere on $\RO{0,t}\times\R_+^{n-1}$ and not just almost everywhere. (Doing so, we might loose some units. But, since we still will obtain a generating subset of units, it does not matter.) We will also assume that $\xi^\odot$ is an exponential unit, i.e.\ $\xi_t^0=\om$.

Fixing $(t_n,\ldots,t_1)\in\RO{0,s+t}\times\R_+^{n-1}$, we show that $\s_t\xi_s^k\odot\xi_t^{n-k}(t_n,\ldots,t_1)\ne0$ for at most one $k=1,\ldots,n$. (Of course, $k$ depends on $(t_n,\ldots,t_1)$.) To see this choose a $k_0$ with nonzero contribution to \eqref{nunit}. Then $t_{n-k_0}$ must be in $\RO{0,t}$. Therefore, $t_{n-k_0}-t<0$ so that certainly the contribution of all terms in \eqref{nunit} with $k>k_0$ is $0$. From this we conclude that for $k\ne k'$ the contribution for $k$ or that for $k'$ must vanish.

Suppose $t_n<t$. Then only $k=0$ can contribute to \eqref{nunit}. We conclude $\xi_{s+t}=\xi_t^n$ on $\RO{0,t}\times\R_+^{n-1}$. It follows that $\xi_t^n(t_n,\ldots,t_1)$ does not depend on $t$ as long as $t>t_n$ so that there exists a well-defined function $\xi^n\colon\R_+^n\rightarrow F^{\sbar{\odot}n}$ such that $\xi_t^n(t_n,\ldots,t_1)=\I_{\RO{0,t}}(t_n)\xi^n(t_n,\ldots,t_1)$.

Next suppose that $t_\ell\ge t_n\ge t$ for all $\ell$. Then only $k=n$ can contribute to \eqref{nunit}. Inserting $\xi^n$ for $\xi_{s+t}^n$ and $\xi_s^n$ and, finally, putting $t=t_n$, we find $\xi^n(t_n,\ldots,t_1)=\xi^n(0,t_{n-1}-t_n,\ldots,t_1-t_n)$. In other words, there exists a function $\zeta^n\colon\R_+^{n-1}\rightarrow F^{\sbar{\odot}n}$ such that
\beq{ \label{spectup}
\xi_t^n(t_n,\ldots,t_1)
~=~
\I_{\RO{0,t}}(t_n)\zeta^n(t_{n-1}-t_n,\ldots,t_1-t_n)
}\eeq
whenever $t_\ell\ge t_n$ for all $\ell$.

It is clear that defining $\xi_t^n$ on tuples $t_\ell\ge t_n$ as in \eqref{spectup} for an arbitrary family of functions $\zeta^n$ $(n\in\N)$, every part of \eqref{nunit} where $\xi^n$ appears on the right-hand side is satisfied. For tuples not fulfilling $t_\ell\ge t_n$, \eqref{nunit} becomes a recursion to reduce the definition of $\xi_t^n$ to that of $\xi_t^k$ $(k<n)$. In order to obtain a complete definition of $\xi_t$ in terms of all $\zeta^n$, we must decompose an arbitrary tuple $(t_n,\ldots,t_1)\in\RO{0,t}\times\R_+^{n-1}$ into subtuples fulfilling the condition such that \eqref{spectup} can be applied. In order to have well-definedness the decomposition must be unique. The following proposition settles both problems and justifies the \it{Ansatz} in the consecutive theorem.

\bprop
Let  $(t_n,\ldots,t_1)\in\R_+^n$ $(n\in\N)$. Then there exist unique $m,k_1,\ldots,k_m\in\N$ and $s_k^\ell\in\R_+$ $(1\le\ell\le m,1\le k\le k_\ell)$, fulfilling $s_k^\ell\ge s_{k_\ell}^\ell$ $(1\le k\le k_\ell)$, $s_{k_\ell}^\ell>s_{k_{\ell-1}}^{\ell-1}$ and
\beq{ \label{tuple}
(t_n,\ldots,t_1)
~=~
(s_{k_m}^m,\ldots,s_1^m,~\ldots~,s_{k_1}^1,\ldots,s_1^1).
}\eeq
Obviously, every tuple on the right-hand side with $s_i^j$ fulfilling the stated conditions may appear.
\eprop

\proof
For $n=1$ the statement is clear. We proceed by induction on $n$. Start with removing $t_n$ from the tuple and, going backwards, continue removing all $t_k$, as long as $t_k\ge t_n$. If the remaining tuple is empty, then we are done. If the remaining tuple is nonempty, then its left entry is strictly smaller than $t_n$. Applying the induction hypothesis to this tuple (whose length is smaller than $n$), existence follows.

Once existence is established, it follows that the first tuple to be removed from the left is unique. Hence, uniqueness follows, once more, by induction.~\qed

\lf
The points of a typical tuple as in \eqref{tuple} may be visualized as in the following diagram. In each subtuple $(s_{k_\ell}^\ell,\ldots,s_1^\ell)$ the first element $s_{k_\ell}^\ell$ must hit the thick line while the remaining ones must be (not necessarily strictly) above the thick line.

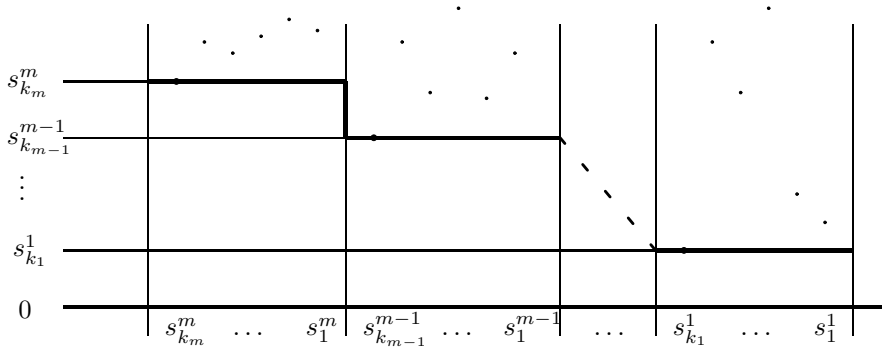
\begin{figure}[H]
\centering
\begin{picture}(155,60)(-10,-5)
\Thicklines

\put(0,0){\line(1,0){145}}
\put(-8,-2){$0$}
\put(15,40){\line(1,0){35}}
\put(50,40){\line(0,-1){10}}
\put(50,30){\line(1,0){38}}
\put(105,10){\line(1,0){35}}
\Thicklines
\dashline{2}(88,30)(105,10)

\thinlines

\put(15,-5){\line(0,1){55}}

\put(0,40){\line(1,0){15}}
\put(-10,39){$s_{k_m}^m$}
\put(20,40){\circle*{1}}
\put(18,-5){$s_{k_m}^m$}
\put(25,47){\circle*{.5}}
\put(30,45){\circle*{.5}}
\put(35,48){\circle*{.5}}
\put(40,51){\circle*{.5}}
\put(45,49){\circle*{.5}}
\put(43,-5){$s_1^m$}
\put(30,-5){$\ldots$}

\put(50,-5){\line(0,1){55}}

\put(0,30){\line(1,0){50}}
\put(-10,29){$s_{k_{m-1}}^{m-1}$}
\put(55,30){\circle*{1}}
\put(53,-5){$s_{k_{m-1}}^{m-1}$}
\put(60,47){\circle*{.5}}
\put(65,38){\circle*{.5}}
\put(70,53){\circle*{.5}}
\put(75,37){\circle*{.5}}
\put(80,45){\circle*{.5}}
\put(78,-5){$s_1^{m-1}$}
\put(67,-5){$\ldots$}

\put(88,-5){\line(0,1){55}}

\put(-8,19){$\vdots$}
\put(94,-5){$\ldots$}

\put(105,-5){\line(0,1){55}}

\put(0,10){\line(1,0){105}}
\put(-9,9){$s_{k_1}^1$}
\put(110,10){\circle*{1}}
\put(108,-5){$s_{k_1}^1$}
\put(115,47){\circle*{.5}}
\put(120,38){\circle*{.5}}
\put(125,53){\circle*{.5}}
\put(130,20){\circle*{.5}}
\put(135,15){\circle*{.5}}
\put(133,-5){$s_1^1$}
\put(120,-5){$\ldots$}

\put(140,-5){\line(0,1){55}}

\end{picture}
\end{figure}

\bcor \label{arrcor}
Let $S_t$ be the disjoint union of all $\RO{0,t}\times\R_+^{n-1}$ $(n\in\N)$ and let $f=\bfam{f_n}_{n\in\N}$ be a summable family of integrable functions $f_n\colon\RO{0,t}\times\R_+^{n-1}\rightarrow\cB$. (That means $f\in L^1(S_t,\cB)$.) Then
\bmun{
\int_{S_t}f\,dS
~:=~
\sum_{n=1}^\infty\int_0^t\,dt_n\,\int_0^\infty\,dt_{n-1}\,\ldots\int_0^\infty\,dt_1\,f_n(t_n,\ldots,t_1)
\\
~=~
\sum_{n=1}^\infty~\sum_{k_1+\ldots+k_m=n}\int_0^t\,ds_{k_m}^m\,\int_{s_{k_m}}^\infty\,ds_{k_m-1}^m\,\ldots\int_{s_{k_m}^m}^\infty\,ds_1^m\,~~~~~~~~~~~~~~~~~~~
\\
\ldots~\int_0^{s_{k_2}^2}\,ds_{k_1}^1\,\int_{s_{k_1}}^\infty\,ds_{k_1-1}^1\,\ldots\int_{s_{k_1}^1}^\infty\,ds_1^1\,f_n(s_{k_m}^m,\ldots,s_1^m,~\ldots~,s_{k_1}^1,\ldots,s_1^1).
}\emun
After a resummation and substitutions $t_\ell$ for $s_{k_\ell}^\ell$ and $s_k^\ell$ for $s_k^\ell-t_\ell$ $(k<k_\ell)$, we obtain
\bmun{
\sum_{m=1}^\infty~\sum_{k_1,\ldots,k_m=1}^\infty\int_0^t\,dt_m\,\int_0^{t_m}\,dt_{m-1}\,\ldots~\int_0^{t_2}\,dt_1\,
\\
\int_0^\infty\,ds_{k_m-1}^m\,\ldots\int_0^\infty\,ds_1^m\,~\ldots~\int_0^\infty\,ds_{k_1-1}^1\,\ldots\int_0^\infty\,ds_1^1\,~~~~~~~~~~~~~~~~
\\
f_{k_1+\ldots+k_m}(t_m,s_{k_m-1}^m+t_m,\ldots,s_1^m+t_m,~\ldots~,t_1,s_{k_1-1}^1+t_1,\ldots,s_1^1+t_1).
}\emun
\ecor

\bthm\label{Fowmodthm}
Let $\zeta=\bigoplus_{n\in\N}\zeta^n\in\coplus_{n\in\N}L^2(\R_+^{n-1},F^{\sbar{\odot}n})$ $(=F\sbar{\odot}\cF(E_{\R_+})=\cF(E_{\R_+})\sbar{\odot}F)$. Then
\beqn{
\xi_t^n(t_n,\ldots,t_1)
~=~
\zeta^{k_m}(s_{k_m-1}^m-s_{k_m}^m,\ldots,s_1^m-s_{k_m}^m)
\odot\ldots\odot
\zeta^{k_1}(s_{k_1-1}^1-s_{k_1}^1,\ldots,s_1^1-s_{k_1}^1)
}\eeqn
defines an exponential unit ${\xi^\zeta}^\odot$ with $\xi^\zeta_t=\bigoplus_{n\in\N}\xi_t^n$. Moreover,  the subset $\bCB{{\xi^\zeta}^\odot}$ of units contains all exponential units of $\cF^\odot(F)$, it is generating (so that $\cF^\odot(F)$ is completely spatial) and ${\xi^\zeta}^\odot\mapsto\xi^\odot(0,\zeta)$ defines an isomorphism $\cF^\odot(F)\rightarrow\DG^\odot(F\sbar{\odot}\cF(E_{\R_+}))$. In other words, $\cF^\odot(F)$ is that unique completely spatial product system with index $F\sbar{\odot}\cF(E_{\R_+})$
\ethm

\proof
First, we show that, if all $\xi^\zeta_t$ are in $\cF_t(F)$, then they define a unit. For that it is sufficient to show that they fulfill \eqref{nunit}. Let us first consider only the case $m=1$. This is precisely the case where, depending on whether $t_n<t$ or $t_n\ge t$, only $k=0$ or $k=n$ contribute. In both cases equality of left- and right-hand side is immediate. Taking into account that $\cF^\odot(F)$ is a product system, a short consideration shows that the general case follows similarly, but requires harder work in order to write it down. We omit this.

Next, we show that ${\xi^\zeta}^\odot\mapsto\xi^\odot(0,\zeta)$ is isometric and \nbd{\cB}\nbd{\cB}linear. This establishes, in particular, that $\xi_t$ is an element of $\cF_t(F)$. We have to show that given $\zeta,\zeta'$ and $b$, we have $\AB{\xi^\zeta_t,b\xi^{\zeta'}_t}=\AB{\xi_t(0,\zeta),b\xi(0,\zeta')}$.
\beqn{
\AB{\xi^\zeta_t,b\xi^{\zeta'}_t}
~=~
\sum_{n=0}^\infty\int_0^t\,dt_n\,\int_0^\infty\,dt_{n-1}\,\ldots\int_0^\infty\,dt_1\,\AB{\xi^n_t,b{\xi'_t}^n}(t_n,\ldots,t_1).
}\eeqn
Inserting the concrete form of $\xi_t^n$, by Corollary \ref{arrcor} we obtain
\bmun{
\U+\sum_{m=1}^\infty~\sum_{k_1,\ldots,k_m=1}^\infty\int_0^t\,dt_m\,\int_0^{t_m}\,dt_{m-1}\,\ldots~\int_0^{t_2}\,dt_1\,
\\
\int_0^\infty\,ds_{k_m-1}^m\,\ldots\int_0^\infty\,ds_1^m\,~\ldots~\int_0^\infty\,ds_{k_1-1}^1\,\ldots\int_0^\infty\,ds_1^1\,~~~~~~~~~~~~
\\
\BAB{\zeta^{k_1}~\ldots~,~\BAB{\zeta^{k_m},b{\zeta'}^{k_m}}(s_{k_m-1}^m,\ldots,s_1^m)~\ldots~{\zeta'}^{k_1}}(s_{k_1-1}^1,\ldots,s_1^1)
\\[2ex]
~=~
\U+\sum_{m=1}^\infty\frac{t^m}{m!}~\sum_{k_1,\ldots,k_m=1}^\infty\AB{\zeta^{k_m}\odot\ldots\odot\zeta^{k_1},b{\zeta'}^{k_m}\odot\ldots\odot{\zeta'}^{k_1}}
\\[1ex]
~=~
\U+\sum_{m=1}^\infty\frac{t^m}{m!}\AB{\zeta^{\odot m},b{\zeta'}^{\odot m}}
~=~
\AB{\xi_t(0,\zeta),b\xi_t(0,\zeta')}.
}\emun
Therefore, ${\xi^\zeta}^\odot\mapsto\xi^\odot(0,\zeta)$ is an isomorphism from the subsystem generated by all ${\xi^\zeta}^\odot$ onto $\DG^\odot(F\sbar{\odot}\cF(E_{\R_+}))$.

What remains to show is that $\CB{{\xi^\zeta}^\odot}$ is generating, for in this case $\cF^\odot(F)$ is isomorphic to $\DG^\odot(F\sbar{\odot}\cF(E_{\R_+}))$ and, in particular,  the former cannot have more exponential units, because the latter does not have more.

Recall the well-known fact (see, e.g., [\refcite{Ske01}, Chapter 7]) that the functions
\beqn{
x_n
\colon
(t_n,\ldots,t_1)
~\longmapsto~
\I_{\RO{r_n,s_n}}(t_n)\ldots\I_{\RO{r_1,s_1}}(t_1)y_n\odot\ldots\odot y_1
}\eeqn
$(r_\ell<s_\ell,s_n<t,y_\ell\in F)$ with $\RO{r_\ell,s_\ell}\cap\RO{r_{\ell'},s_{\ell'}}=\emptyset$ for $\ell\ne\ell'$, form a total subset of $L^2(\RO{0,t}\times\R_+^{n-1},F^{\sbar{\odot}n})$. Such a function is supported by a set of tuples for which the parameters $m,k_1,\ldots,$ $k_m\in\N$ in \eqref{tuple} are fixed. Moreover, it decomposes into a tensor product of functions with $m=1$. Since the tensor product is linear in its factors, it is sufficient to understand only the case $m=1,k_1=n$, i.e.\ $t_n\in\RO{0,t},t_k\ge t_n$ for all $k$. We may even assume that $\RO{r_n,s_n}=\RO{0,t}$. (Otherwise, factorize into a tensor product with $\om_{t-s_n}$ and $\om_{r_n}$.) But in this case we find easily by differentiating $\xi^{\lambda\zeta}_t$ with respect to $\lambda$ at $\lambda=0$ for $\zeta(t_{n-1},\ldots,t_1)=x_n(t_n,\ldots,t_1)$ that $x_n$ is contained in what $\CB{{\xi^\zeta}^\odot}$ generates.~\qed

\lf
Choosing for $F$ a (separable) Hilbert space $K\ne\CB{0}$, we recover Fowler's result that the index of $\sS$ is $\dim(K\sbar{\otimes}\cF(E_{\R_+}))=\infty$. In Fowler's case the ``number'' $\infty$ was just sufficient, and although also in [\refcite{Fow95}] the isomorphism was constructed explicitly (very much like we do here, but somewhat much more complicated), in the formulation of the main result the concrete structure of the index space was neglected. In the module case the index is the whole space and we cannot stick to a simple dimension. The index of a free flow depends, in general, on its free index.

\section{Open problems and outlook} \label{openSEC}

{\parskip0.5ex plus 0.5ex minus 0.5ex
\bemp[Quantum stochastic calculus on spatial product systems?~]\label{calcop}
Consider the triples $(E,\vt,\xi)$ and $(E',\vt',\xi')$ (not necessarily dilations or noises). If $E\cong E'$, then by [\refcite{Ske01}, Theorem 14.1.5] the triples have isomorphic product systems, if and only the \nbd{E_0}semigroups $\vt$ and $\vt'$ are \it{cocycle conjugate}. By [\refcite{Ske01}, Proposition 14.1.2] this is independent of the choice of the unit vectors $\xi,\xi'$ and also the version without unit vector based on [\refcite{MSS03p}] admits such a uniqueness result as long $E\cong E'$ are full. (In the case of infinite-dimensional separable Hilbert spaces, as considered by Arveson, $E$ and $E'$ are always isomorphic.) It is not difficult to see that isomorphic product systems need not imply that $E\cong E'$, but examples that can be constructed easily are not primary. In the case of primary noises it is rather easy to show that $E\cong E'$, if there is a spatial isomorphism of the product systems. (Loosely speaking, the product systems coincide and there is an automorphism that sends one reference unit to the other. Once again, we meat the question whether the automorphism group of a product system acts transitively on the central units. If the product system is completely spatial, then we know from Remark \ref{transcsrem} that this is true. In cases when a weak dilation of a CP-semigroup is constructed with the help of a quantum stochastic calculus on a symmetric Fock module [\refcite{Bha01}] or a full Fock module [\refcite{Ske00}] considered as inductive limit over the spatial vacuum unit, we know that the inductive limit with respect to the unit $\xi^\odot$ that generates the CP-semigroup as $T_t=\AB{\xi_t,\bullet\xi_t}$, a noncentral unit if $T$ is nontrivial, coincides or is contained in the Fock module. It would be, indeed, a great achievment to show in full generality that, if the product system $E^\odot$ with the unit $\xi^\odot$ of a possibly nontrivial CP-semigroup $T$ is spatial, with reference unit $\om^\odot$ say, then the primary weak dilation of $T$ contructed as inductive limit over the unit $\xi^\odot$ is a cocycle perturbation of the primary noise constructed as inductive limit over the reference unit $\om^\odot$. Presently, we are only dreaming of a calculus on spatial product sytems, permitting to prove such a result.) One may check for two pairs $(E^\ell,\vt^\ell,\om^\ell)$ and $(E^\ell,{\vt'}^\ell,{\om'}^\ell)$ $(\ell=1,2)$ of cocycle conjugate noises that also the product noises $(\bfam{E^1\circledcirc E^2},\bfam{\vt^1\circledcirc\vt^2},\om)$ and $(\bfam{E^1\circledcirc'E^2}\bfam{{\vt'}^1\circledcirc{\vt'}^2},\om')$ are cocycle conjugate. In particular, the product systems $\bfam{E^1\circledcirc E^2}^{\sbar{\odot}}$ and $\bfam{E^1\circledcirc'E^2}^{\sbar{\odot}}$ constructed from the reference units ${\om^\ell}^\odot$ and ${{\om'}^\ell}^\odot$, respectively, are isomorphic (but not necessarily as spatial systems).
\eemp

\bemp[Unital embeddings, unital dilations and free product systems?~]\label{freeps}
In the case of Hilbert spaces the tensor product  of two \nbd{E_0}semigroups $(\eH^\ell,\vt^\ell)$ is the $E_0$--semi\-group $(\eH^1\sbar{\otimes}\eH^2,\vt^1\otimes\vt^2)$. The algebra $\sB(\eH^1\sbar{\otimes}\eH^2)$ contains the original algebras $\sB(\eH^\ell)$ as unital subalgebras $\sB(\eH^1)\otimes\id_{\eH^2}$ and $\id_{\eH^1}\otimes\sB(\eH^2)$, and the restrictions of $\vt^1\otimes\vt^2$ to these subalgebras gives back $\vt^\ell$. In the case of Hilbert modules we do not know unital embeddings $\sB^a(E^\ell)\rightarrow\sB^a(\bfam{E^1\circledcirc E^2})$. Of course, we may embed $\sB^a(E^\ell)$ nonunitally with the help of the projection $p^\ell$ onto the submodule $E^\ell$ of $\bfam{E^1\circledcirc E^2}$. However, the projections $\bfam{\vt^1\circledcirc\vt^2}_t(p^\ell)$ are, in general, strictly increasing so that $\bfam{\vt^1\circledcirc\vt^2}_t$ does not leave invariant the embedded unit $p^\ell$ of $\sB^a(E^\ell)$. (Only the compression $p^\ell\bfam{\vt^1\circledcirc\vt^2}p^\ell\upharpoonright p^\ell\sB^a(\bfam{E^1\circledcirc E^2})p^\ell$ gives back $\vt^\ell$. Notice also that the two nonunital subalgebras $\sB^a(E^1)$ and $\sB^a(E^2)$ are boolean indipendent in the conditional expectation $\AB{\om,\bullet\om}$ in the sense of Skeide [\refcite{Ske00}, Section 14].)

A possible solution for a product of noises with unital embeddings is suggested by a comparison of CCR-flows and free flows as discussed in Section \ref{freeflow}. The free flow with free index $F^1\oplus F^2$ is nothing but the free product $\sS^1*\sS^2$ of the free flows $\sS^\ell$ with index $F^\ell$ in the vacuum conditional expectation in the sense of Voiculescu [\refcite{Voi95}] and Speicher [\refcite{Spe98}]. Clearly, there exist unital embeddings of $\sB^a(\cF(E^\ell_{\R_+}))$ into $\sB^a(\cF(E^1_{\R_+}\oplus E^2_{\R_+}))$ which behave covariantly under the relevant time shifts. Observe that $\cF(E_{\RO{0,s}})*\cF(E_{\RO{0,t}})=\cF(E_{\RO{0,s+t}})$ where $*$ indicates the free product of Hilbert \nbd{\cB}\nbd{\cB}modules with respect to the reference vectors $\om_s,\om_t$. In other words, we may associate a free product system with each free flow and the free product of the members of the free product systems for $\sS^1$ and $\sS^2$ gives us the free product system of the free product $\sS^1*\sS^2$. Free flows are primary noises and similarly as for (tensor) product systems we may recover the flow by an inductive limit over the members of its free product system. It can be shown that every spatial product system generates a universal free product system into which it embeds and that is generated by it as a free product system. The free product system of the full Fock module can be interpreted in this sense as the free product system generated by the time ordered product subsystem. This raises several questions: Is every free product system the free product system generated by a spatial product system? Is there a general possibility to decide whether a primary noise stems from a free product system? (Of course, we know that it stems from a spatial product system.) Most important, and even for Fock modules not yet answered: Can we understand a weak dilation with a spatial product system as a projection from a unital dilation to the free flow associated with the free product system generated by the spatial product system? (To answer this question it will definitely by important to answer first the question about cocycle conjugacy of the weak dilation with spatial product system and the noise associated with that spatial product system as raised in \ref{calcop}.) We will leave such questions to future work.
\eemp

\bemp[An index for nonspatial product systems?~]\label{indnonspat}
What is the index of a nonspatial system? For Arveson systems the answer is motivated by the analogue of Theorem \ref{uprod} which asserts that units in the product may be composed from units in the factors and remains true for arbitrary Arveson systems. Therefore, if one factor has no units, then also the tensor product has no units. Putting the index of such systems to $\infty$ (better: $-\infty$, because in this case the index has somewhat from $\log0$), the computation rules for the Arveson index (in the spatial case the dimension of $K$ for the maximal completely spatial subsystem $\DG^\otimes(K)$) remain valid. Thus, we could introduce a fomal space $\text{``$\infty$''}$ as index of a nonspatial system with the formal computation rule $\text{``$\infty$''}\oplus F=\text{``$\infty$''}$. However, as long as we do not have an extension of the product $\circledcirc$ to nonspatial subsystems, such a definition is meaningless. A way out could be to transform a nonspatial system into a spatial one. For instance, by the first attempt of defining a product, as described in the beginning of Section \ref{prodSEC}, where existing units in the factors are mapped to orthogonal units in the product, we could construct this product for a nonspatial product system and the trivial one. This product contains a central unital unit, namely, the unit $\bfam{\U}_{t\in\R_+}$ of the trivial product system, and, therefore, is spatial. We may construct the product $\circledcirc$ of this extended product system with any other spatial one, and then restrict to what the original product system generates there. Also this product system should be nonspatial, because otherwise by Theorem \ref{uprod} there must be a unit also in the original nonspatial factor. We postpone a detailed analysis of this idea. But, it is clear that it cannot serve as a general construction, because if we apply it to systems which are already spatial, then two units from different factors will be orthogonal. After all, it might appear not to be very reasonable to extend the definition of index to nonspatial systems. (Already in the case of Arveson systems the index $\infty$ of type III systems does not mean that there is a maximal completely spatial subsystem $\DG^\otimes(K)$ for some infinite-dimensional Hilbert space and, therefore, the concrete meaning of the index is obstructed by that definition.)
\eemp

\bemp[An example of Powers naturally leading to our product of spatial product systems.~]\label{powersprod}
We close with the remark that recently an interesting interpretation of our product has occured. In the 2002 AMS-Meeting ``Advances in Quantum Dynamics'' R.\ Powers constructed a CP-semigroup on $\sB(H\oplus H)$ from two \nbd{E_0}semigroups $\vt^i$ on $\sB(H)$ with spatial product systems ${\eH^i}^\otimes$ and reference units ${\om^i}^\otimes$ by setting
\beqn{
T_t
\Matrix{
a_{11}&a_{12}
\\
a_{21}&a_{22}
}
~=~
\Matrix{
a_{11}\otimes\id_{\eH^1_t}&(\id_H\otimes\om^1_t)a_{12}(\id_H\otimes\om^2_t)^*
\\
(\id_H\otimes\om^2_t)a_{21}(\id_H\otimes\om^1_t)^*&a_{22}\otimes\id_{\eH^2_t}
}.
}\eeqn
(We make use of the identifications $H\otimes\eH^1_t=H=H\otimes\eH^2_t$. For instance, $a_{ii}\otimes\id_{\eH^i_t}$ is just $\vt^i_t(a_{ii})$.) Powers asked for the product system of that CP-semigroup. (The product system of a CP-semigroup was defined by Bhat [\refcite{Bha96}] as product system of the minimal dilating \nbd{E_0}semigroup. With the methods from Bhat and Skeide [\refcite{BhSk00}] it is possible to contruct the product system of a CP-semigroup directly without dilation.) Still during that meeting we were able to show (see Skeide [\refcite{Ske03c}]) that the product system of that CP-semigroup is exactly our product of the spatial product systems of the two \nbd{E_0}semigroups with respect to the reference units. By Remark \ref{ntprem} it may be but it need not be the tensor product (answering Powers' question about the product system posed at the meeting). In Bhat, Liebscher and Skeide [\refcite{BLS02p}] we will generalize Powers' example and the new interpretation in [\refcite{Ske03c}] from $\sB(H)$ to $\sB^a(E)$.
\eemp
}

\section*{Appendix: Lemma \ref{Tlem} and variations}\label{APP}

\renewcommand{\thesection}{A}
\setcounter{equation}{0}

In this appendix we proof the following variation for net limits over $\bJ_t$ Lemma \ref{Tlem} for net limits over $\bJ_t$ but only in the case of spatial product systems.\footnote{Meanwhile, a full proof without reference to time ordered product systems has appeared in Liebscher and Skeide [\refcite{LiSk05p}].}

\blem \label{netlem}
Let ${\xi^\ell}^\odot$ $(\ell=1,2)$ be units in a continuous subset $S$ of units in a product system $E^{\sbar{\odot}}$. Then for all $\vk^1,\vk^2\in\C$ with $\vk^1+\vk^2=1$ the limit
\beq{\label{netTrotter}
\xi_t
~=~
\lim_{\et\in\bJ_t}(\vk^1\xi^1_{t_n}+\vk^2\xi^2_{t_n})\odot\ldots\odot(\vk^1\xi^1_{t_1}+\vk^2\xi^2_{t_1})
}\eeq
exists in norm, $\xi^\odot=\bfam{\xi_t}_{t\in\R_+}$ is a unit, too, and the set $S\cup\CB{\xi^\odot}$ is still continuous. Moreover, for all ${\xi'}^\odot\in S\cup\CB{\xi^\odot}$ we have $\eL^{\xi',\xi}=\vk^1\eL^{\xi',\xi^1}+\vk^2\eL^{\xi',\xi^2}$.
\elem

The limit turns out to be uniform in the norm $\norm{\et}=\max\CB{t_n,\ldots,t_1}$ of $\et\in\bJ_t$. This proves, in particular, also Lemma \ref{Tlem}. Here we will proof only the spatial case, where all units are known so that we have a concrete candidate for the limit.

\proof
By Theorem \ref{msssthm} we may assume that we are speaking about units ${\xi^\ell}^\odot=\xi^\odot(\beta_\ell,\zeta_\ell)$ in a time ordered product system. We shall denote $y_t:=\vk^1\xi_t(\beta_1,\zeta_1)+\vk^2\xi_t(\beta_2,\zeta_2)$ and $y_\et:=y_{t_n}\odot\ldots\odot y_{t_1}$. A candidate for the limit is the unit $\xi^\odot(\vk^1\beta_1+\vk^2\beta_2,\vk^1\zeta_1+\vk^2\zeta_2)$ whose elements we denote by $z_t$. If we show convergence to that limit, then the remaining statements follow. We will show
\beq{
\label{L1}
\AB{y_\et,y_\et}
~\longrightarrow~
\AB{z_t,z_t}
}\eeq
and
\beq{
\label{L2}
\AB{\xi_t(\beta,\zeta),y_\et}
~\longrightarrow~
\AB{\xi_t(\beta,\zeta),bz_t}
\text{~~~for all~~~}
(\beta,\zeta),b.
}\eeq
\eqref{L1} implies, in particular, that the net $y_\et$ is eventually bounded. (Actually, the net is bounded, what we must show first, and all convergences are in the norm of mappings $\AB{x,\bullet y}$ in $\sB(\cB)$.) As elements like in \eqref{tensgen} are total, (enventual) boundedness of the $y_\et$ implies $\AB{x,y_\et}\to\AB{x,z_t}$ for all $x\in E_t$. Therefore, like for Hilbert spaces, it follows
\beqn{
\AB{z_t-y_\et,z_t-y_\et}
~=~
\AB{z_t,z_t}+\AB{y_\et,y_\et}-\AB{z_t,y_\et}-\AB{y_\et,z_t}
~\longrightarrow~
0.
}\eeqn

We compute
\baln{
\AB{y_t,\bullet y_t}
~=~
&
\bar{\vk}_1\vk^1\eU_t^{(\beta_1,\zeta_1)(\beta_1,\zeta_1)}+
\bar{\vk}_1\vk^2\eU_t^{(\beta_1,\zeta_1)(\beta_2,\zeta_2)}
\\
&
~~~~~~~~~~~~+
\bar{\vk}_2\vk^1\eU_t^{(\beta_2,\zeta_2)(\beta_1,\zeta_1)}+
\bar{\vk}_2\vk^2\eU_t^{(\beta_2,\zeta_2)(\beta_2,\zeta_2)}
\\
~=~
&
(\bar{\vk}_1\vk^1+\bar{\vk}_1\vk^2+\bar{\vk}_2\vk^1+\bar{\vk}_2\vk^2)\id_\cB
\\
&
~~~~~~+
t\Bigl(\bar{\vk}_1\vk^1\eL^{(\beta_1,\zeta_1)(\beta_1,\zeta_1)}+
\bar{\vk}_1\vk^2\eL^{(\beta_1,\zeta_1)(\beta_2,\zeta_2)}
\\
&
~~~~~~~~~~~~+
\bar{\vk}_2\vk^1\eL^{(\beta_2,\zeta_2)(\beta_1,\zeta_1)}+
\bar{\vk}_2\vk^2\eL^{(\beta_2,\zeta_2)(\beta_2,\zeta_2)}\Bigr)
+
O(t^2)
\\
~=~
&
\id_\cB+t\eL^{(\vk^1\beta_1+\vk^2\beta_2,\vk^1\zeta_1+\vk^2\zeta_2),(\vk^1\beta_1+\vk^2\beta_2,\vk^1\zeta_1+\vk^2\zeta_2)}+O(t^2)
\\
~=~
&
\AB{z_t,\bullet z_t}+O(t^2)
}\emun
with an $O(t^2)$ that may be chosen uniformly in $t$ on every compact interval $\SB{0,T}$. Let us denote $Y_t=\AB{y_t,\bullet y_t}$ and $Z_t=\AB{z_t,\bullet z_t}$. By standard arguments there is a global constant $c>0$ such that both $Y_{t_1}\circ\ldots\circ Y_{t_n}$ and $Z_{t_1}\circ\ldots\circ Z_{t_n}$ can be estimated by $e^{ct}$ for every $t\ge0$ and every $\et\in\bJ_t$. We find
\baln{
\AB{y_\et,\bullet y_\et}-\AB{z_t,\bullet z_t}
~=~
&
Y_{t_1}\circ\ldots\circ Y_{t_n}-Z_{t_1}\circ\ldots\circ Z_{t_n}
\\
~=~
&
\sum_{k=1}^nY_{t_1}\circ\ldots\circ Y_{t_{k-1}}\circ(Y_{t_k}-Z_{t_k})\circ Z_{t_{k+1}}\circ\ldots\circ Z_{t_n}
}\ealn
whose norm can be estimated by
\beqn{
M_t\sum_{k=1}^nt_k^2
~\le~
M_t\norm{\et}\sum_{k=1}^nt_k
~=~
M_t\norm{\et}t.
}\eeqn
This shows \eqref{L1}. By completely analogue but simpler computations we show also \eqref{L2}.~\qed

\newcommand{\Swap}[2]{#2#1}\newcommand{\Sort}[1]{}
\providecommand{\bysame}{\leavevmode\hbox to3em{\hrulefill}\thinspace}
\providecommand{\MR}{\relax\ifhmode\unskip\space\fi MR }
\providecommand{\MRhref}[2]{%
  \href{http://www.ams.org/mathscinet-getitem?mr=#1}{#2}
}
\providecommand{\href}[2]{#2}





\begin{thebibliography}{10}

\bibitem{AcKo01}
L.~Accardi and S.~Kozyrev, \emph{{On the structure of Markov flows}}, Chaos
  Solitons Fractals \textbf{12} (2001), 2639--2655.

\bibitem{Arv89}
W.~Arveson, \emph{{Continuous analogues of Fock space}}, Mem. Amer. Math. Soc.,
  no. 409, American Mathematical Society, 1989.

\bibitem{Arv89a}
\bysame, \emph{{Continuous analogues of Fock space III: Singular states}}, J.\
  Operator Theory \textbf{22} (1989), 165--205.

\bibitem{BBLS04}
S.D. Barreto, B.V.R. Bhat, V.~Liebscher, and M.~Skeide, \emph{{Type I product
  systems of Hilbert modules}}, J.\ Funct.\ Anal. \textbf{212} (2004),
  121--181.

\bibitem{Bha96}
B.V.R. Bhat, \emph{{An index theory for quantum dynamical semigroups}}, Trans.\
  Amer.\ Math.\ Soc. \textbf{348} (1996), 561--583.

\bibitem{Bha01}
\bysame, \emph{{Cocycles of CCR-flows}}, Mem. Amer. Math. Soc., no. 709,
  American Mathematical Society, 2001.

\bibitem{BLS02p}
B.V.R. Bhat, V.~Liebscher, and M.~Skeide, \emph{{(Tentative title) A problem of
  Powers and the product of spatial product systems}}, Preprint, Campobasso, in
  preparation, 2004.

\bibitem{BhSk00}
B.V.R. Bhat and M.~Skeide, \emph{{Tensor product systems of Hilbert modules and
  dilations of completely positive semigroups}}, Infin. Dimens. Anal. Quantum
  Probab. Relat. Top. \textbf{3} (2000), 519--575.

\bibitem{BhSr05}
B.V.R. Bhat and R.~Srinivasan, \emph{{On product systems arising from sum
  systems}}, Infin. Dimens. Anal. Quantum Probab. Relat. Top. \textbf{8}
  (2005), 1--31.

\bibitem{ChrEv79}
E.~Christensen and D.E. Evans, \emph{{Cohomology of operator algebras and
  quantum dynamical semigroups}}, J.\ London Math.\ Soc. \textbf{20} (1979),
  358--368.

\bibitem{Fow95}
N.J. Fowler, \emph{{Free $E_0$--semigroups}}, Can.\ J.\ Math. \textbf{47}
  (1995), 744--785.

\bibitem{HKK04p}
J.~Hellmich, C.~K\"ostler, and B.~K\"ummerer, \emph{{Noncommutative continuous
  Bernoulli sifts}}, Preprint, ArXiv: math.OA/0411565, 2004.

\bibitem{Hir04}
I.~Hirshberg, \emph{{$C^*$--Algebras of Hilbert module product systems}}, J.\
  Reine Angew.\ Math. \textbf{570} (2004), 131--142.

\bibitem{Hir05}
\bysame, \emph{{On the universal property of Pimsner-Toeplitz $C^*$--algebras
  and their continuous analogues}}, J.\ Funct.\ Anal. \textbf{219} (2005),
  21--33.

\bibitem{Lan95}
E.C. Lance, \emph{{Hilbert $C^*$--modules}}, Cambridge University Press, 1995.

\bibitem{Lie00p1}
V.~Liebscher, \emph{{Random sets and invariants for (type II) continuous tensor
  product systems of Hilbert spaces}}, Preprint, ArXiv: math.PR/0306365, 2003.

\bibitem{LiSk01}
V.~Liebscher and M.~Skeide, \emph{{Units for the time ordered Fock module}},
  Infin. Dimens. Anal. Quantum Probab. Relat. Top. \textbf{4} (2001), 545--551.

\bibitem{LiSk05p}
\bysame, \emph{{Contructing new units by composing old ones}}, Preprint,
  Campobasso, 2005.

\bibitem{MSS03p}
P.S. Muhly, M.~Skeide, and B.~Solel, \emph{{Representations of $\sB^a(E)$}},
  Preprint, ArXiv: math.OA/0410607, 2004.

\bibitem{MuSo02}
P.S. Muhly and B.~Solel, \emph{{Quantum Markov processes (correspondences and
  dilations)}}, Int.\ J.\ Math. \textbf{51} (2002), 863--906.

\bibitem{MuSo04}
\bysame, \emph{{Hardy algebras, $W^*$--correspondences and interpolation
  theory}}, Math.\ Ann. \textbf{330} (2004), 353--415.

\bibitem{PaSu98}
K.R. Parthasarathy and V.S. Sunder, \emph{{Exponentials of indicater functions
  are total in the boson Fock space $\G(L^2[0,1])$}}, Quantum Probability
  Communications X (R.L. Hudson and J.M. Lindsay, eds.), World Scientific,
  1998, pp.~281--284.

\bibitem{Pas73}
W.L. Paschke, \emph{{Inner product modules over $B^*$--algebras}}, Trans.\
  Amer.\ Math.\ Soc. \textbf{182} (1973), 443--468.

\bibitem{Pow87}
R.T. Powers, \emph{{A non-spatial continuous semigroup of $*$--endomorphisms of
  $\sB(\eH)$}}, Publ. Res. Inst. Math. Sci. \textbf{23} (1987), 1053--1069.

\bibitem{Pow88}
\bysame, \emph{{An index theory for semigroups $*$--endomorphisms of
  $\eB(\eH)$}}, Can.\ Jour.\ Math. \textbf{40} (1988), 86--114.

\bibitem{Pow03}
\bysame, \emph{{Continuous spatial semigroups of completely positive maps of
  $\eB(\eH)$}}, New York J.\ Math. \textbf{9} (2003), 165--269, Available at
  \newline \tt{\footnotesize http://nyjm.albany.edu:8000/j/2003/9-13.html}.

\bibitem{Ske00a}
M.~Skeide, \emph{{Indicator functions of intervals are totalizing in the
  symmetric Fock space $\Gamma(L^2(\R_+))$}}, Trends in contemporary infinite
  dimensional analysis and quantum probability (L.~Accardi, H.-H. Kuo,
  N.~Obata, K.~Saito, {Si Si}, and L.~Streit, eds.), Natural and Mathematical
  Sciences Series, vol.~3, Istituto Italiano di Cultura (ISEAS), Kyoto, 2000,
  Volume in honour of Takeyuki Hida.

\bibitem{Ske00}
\bysame, \emph{{Quantum stochastic calculus on full Fock modules}}, J.\ Funct.\
  Anal. \textbf{173} (2000), 401--452.

\bibitem{Ske01}
\bysame, \emph{{Hilbert modules and applications in quantum probability}},
  Ha\-bi\-li\-ta\-tions\-schrift, Cottbus, 2001, Available at \newline
  \tt{\footnotesize
  http://www.math.tu-cottbus.de/INSTITUT/lswas/\_skeide.html}.

\bibitem{Ske02}
\bysame, \emph{{Dilations, product systems and weak dilations}}, Math.\ Notes
  \textbf{71} (2002), 914--923.

\bibitem{Ske03c}
\bysame, \emph{{Commutants of von Neumann modules, representations of
  $\sB^a(E)$ and other topics related to product systems of Hilbert modules}},
  Advances in quantum dynamics (G.L. Price, B~.M. Baker, P.E.T. Jorgensen, and
  P.S. Muhly, eds.), Contemporary Mathematics, no. 335, American Mathematical
  Society, 2003, pp.~253--262.

\bibitem{Ske03b}
\bysame, \emph{{Dilation theory and continuous tensor product systems of
  Hilbert modules}}, QP-PQ: Quantum Probability and White Noise Analysis XV
  (W.~Freudenberg, ed.), World Scientific, 2003.

\bibitem{Ske04}
\bysame, \emph{{Independence and product systems}}, Recent developments in
  stochastic analysis and related topics (S.~Albeverio, Z.-M. Ma, and
  M.~R\"ockner, eds.), World Scientific, 2004.

\bibitem{Ske03p1}
\bysame, \emph{{Intertwiners, dual quasi orthonormal bases and
  representations}}, in preparation, 2004.

\bibitem{Ske04p}
\bysame, \emph{{Unit vectors, Morita equivalence and endomorphisms}}, Preprint,
  ArXiv: math.OA/0412231, 2004.

\bibitem{Spe98}
R.~Speicher, \emph{{Combinatorial theory of the free product with amalgamation
  and operator-valued free probability theory}}, Mem. Amer. Math. Soc., no.
  627, American Mathematical Society, 1998.

\bibitem{Tsi00p1}
B.~Tsirelson, \emph{{From random sets to continuous tensor products: answers to
  three questions of W.\ Arveson}}, Preprint, ArXiv: math.FA/0001070, 2000.

\bibitem{Tsi00p2}
\bysame, \emph{{From slightly coloured noises to unitless product systems}},
  Preprint, ArXiv: math.FA/0006165, 2000.

\bibitem{Tsi04p1}
\bysame, \emph{{On automorphisms of type II Arveson systems (probabilistic
  approach)}}, Preprint, ArXiv: math.OA/0411062, 2004.

\bibitem{Voi95}
D.~Voiculescu, \emph{{Operations on certain non-commutative operator-valued
  random variables}}, Ast\'erisque \textbf{232} (1995), 243--275.

\end{thebibliography}
\end{document}